\documentclass[a4paper,11pt]{amsart}
\usepackage{amsfonts,amssymb,amsthm,cite,amsmath,amstext}
\usepackage{color}
\usepackage[colorlinks,linkcolor=black,citecolor=black]{hyperref}
\usepackage{comment}
\usepackage{framed}
\usepackage{calligra} 
\usepackage{tikz-cd}

\definecolor{shadecolor}{gray}{0.875}

\newtheorem{thrm}{Theorem}[section]
\newtheorem{thrmx}{Theorem}

\newtheorem{corx}{Corollary}

\newtheorem{lem}[thrm]{Lemma}
\newtheorem{cor}[thrm]{Corollary}
\newtheorem{prop}[thrm]{Proposition}
\newtheorem{conj}[thrm]{Conjecture}
\newtheorem{expec}[thrm]{Expectation}

\theoremstyle{definition}
\newtheorem{defn}[thrm]{Definition}

\newtheorem{exmple}[thrm]{Example}
\newtheorem{rmk}[thrm]{Remark}

\newtheorem{ques}[thrm]{Question}

\newenvironment{claim}
            {\par \bigskip \noindent \textbf{Claim:}}
            {$\Box$ \par \bigskip \noindent}

\DeclareMathOperator{\exc}{Exc}

\DeclareMathOperator{\vol}{vol}
\DeclareMathOperator{\Mov}{Mov}

\DeclareMathOperator{\supp}{supp}

\DeclareMathOperator{\nd}{nd}

\DeclareMathOperator{\Span}{span}

\DeclareMathOperator{\codim}{codim}

\DeclareMathOperator{\image}{Image}

\DeclareMathOperator{\prim}{Prime}
\DeclareMathOperator{\sing}{sing}

\DeclareMathOperator{\Eff}{\overline{Eff}}
\DeclareMathOperator{\Nef}{Nef}

\DeclareMathOperator{\Amp}{Amp}

\DeclareMathOperator{\cone}{cone}

\DeclareMathOperator{\Null}{Null}
\DeclareMathOperator{\eff}{eff}
\DeclareMathOperator{\zar}{zar}

\title[Numerical characterization of the hard Lefschetz classes of dimension two]{Numerical characterization of the hard Lefschetz classes of dimension two, I:
\\
\footnotesize{supercritical collections under certain rearrangement }}
\author{Jiajun Hu and Jian Xiao}
\date{}

\begin{document}
\maketitle


\begin{abstract}
We study the numerical characterization of two dimensional hard Lefschetz classes given by the complete intersections of nef classes. 
In Shenfeld and van Handel's breakthrough work on the characterization of the extremals of the Alexandrov-Fenchel inequality 
for convex polytopes, they proposed an open question on the algebraic analogue of the characterization. 
We settle the open question when the collection of nef classes is given by a rearrangement of supercriticality, 
which in particular includes the big nef collection as a special case.

The main results enable us to refine some previous results and study the extremals of 
Hodge index inequality, and also provide the first series of examples of 
hard Lefschetz classes of dimension two both in algebraic geometry and analytic geometry, 
in which one can allow nontrivial augmented base locus and thus drop the semi-ampleness or semi-positivity assumption.

As a key ingredient of the numerical characterization, we establish a local Hodge index inequality for Lorentzian polynomials, 
which is the algebraic analogue of the local Alexandrov-Fenchel inequality obtained by Shenfeld-van Handel for 
convex polytopes. This result holds in broad contexts, e.g., it holds on a smooth projective variety, 
on a compact K\"ahler manifold, and on a Lorentzian fan, which contains the Bergman fan of a matroid or a polymatroid as a typical example.
\\

\end{abstract}

\tableofcontents

\section{Introduction}

\subsection{Motivation}
As a fundamental result in the topological aspects of an algebraic variety, the hard Lefschetz theorem says that for the cohomology ring of a complex projective manifold of dimension $n$, the $k$-fold product of the cohomology class of a hyperplane $A$ gives an isomorphism:
\begin{equation*}
A^k: H^{n-k}(X) \rightarrow H^{n+k} (X).
\end{equation*}
Since we mostly work with real coefficients thorughout the article, we omit the coefficient here and sometimes in the sequel if not pointed out directly.
It extends to any compact K\"ahler manifolds by replacing the hyperplane class by a K\"ahler class. The classical hard Lefschetz theorem also generalizes vastly to the singular setting and other cohomology theories. In particular, we have the famous Decomposition Theorem \cite{BBDG82DecompsitionThrm}, which in some sense is the deepest result concerning the topological, geometric and arithmetic properties of algebraic varieties, see e.g. \cite{cataldomigDecompSurvey} for a nice survey. The hard Lefschetz theorem is also a fundamental piece of a K\"ahler package. In last decades, when $X$ is a mathematical object other than an algebraic or analytic variety, similar results have been established and playing important roles in algebra, combinatorics and geometry, see e.g. \cite{adiprasito2019combinatorial, huhICM2018, huh2022combinatorics, williamECM, amini2020hodgetropicalArxiv} for the recent developments and the references therein for more classical ones.

The class giving hard Lefschetz type theorems is the main object in this work.
We introduce the following notion: by replacing the $k$-fold product of a hyperplane class or a K\"ahler class by a general $(k,k)$ class $\mathfrak{L}$, we call that $\mathfrak{L}$ is a \emph{hard Lefschetz class} on $H^{n-k}(X)$, if it gives an isomorphism:
\begin{equation*}
\mathfrak{L}: H^{n-k}(X)\rightarrow H^{n+k} (X).
\end{equation*}
We shall utilize the same notion in other cohomology theories if it is applicable, for instance for the groups of numerical equivalent classes of cycles known as the N\'eron-Severi spaces $N^*(X)$. 
Moreover, when Hodge decomposition is available, we call $\mathfrak{L}$ a hard Lefschetz class on $H^{p,q}(X)$ ($p+q=n-k$), if it gives an isomorphism from $H^{p,q}(X)$ to $H^{p+k, q+k} (X)$.
When it does not cause confusion, we just call $\mathfrak{L}$ an \emph{HL class} for short.

In most situations, the HL class $\mathfrak{L}$ arises from strong positivity, such as complete intersections of ample classes, see, for example, \cite{voisinHodge1}, \cite{DN06} and \cite{cattanimixedHRR} for this classical case in algebraic and analytic geometry. For other HL classes $\mathfrak{L}$ that originate from ampleness, we refer the reader to \cite{blochgieseker} and \cite{rosstomaHRRasens, rosstomaSchurLinearImrn}.
It is therefore natural to ask what can be said in the absence of strong positivity assumptions. In previous works \cite{xiaoHodgeIndex, xiaomixedHRR, huxiaohardlef2022Arxiv, hushangxiao2023hardlefArxiv}, several results are obtained towards this question:

\begin{ques}\label{ques HL}
Let $X$ be a smooth projective variety (or a compact K\"ahler manifold) of dimension $n$, and let $\mathfrak{L}$ be a real $(k,k)$-class. When is $\mathfrak{L}$ an HL class?
In a more accessible situation, let $L_1,....,L_{k}$ be real divisor or $(1,1)$-classes. Under what assumptions is the complete intersection $\mathfrak{L}=L_1 \cdot.... \cdot L_{k}$ an HL class?
\end{ques}

One of the main motivations for this question is that, in certain problems and applications, such strong positivity assumptions may not be available or applicable, while hard Lefschetz-type theorems and K\"ahler package still play key roles, see e.g. \cite{adiprasito2019combinatorial, cataldoMigHLsemismallmap, catalMigHodgeMaps}.

In this work, we mainly focus on the case when $\mathfrak{L}$ is a two dimensional class given by the complete intersection of nef classes. This case is closely related to the Hodge index theorem and the log-concavity phenomenon \cite{lazarsfeldPosI, huhICM2018}. Apart from the proposal of Shenfeld and van Handel \cite{handel2022AFextremals} from the perspective of convex geometry and the advances made in \cite{huxiaohardlef2022Arxiv, hushangxiao2023hardlefArxiv} from the viewpoint of algebraic geometry, even a conjectural picture for the characterization of HL classes had been missing. For the general case, we refer the reader to \ref{sec arbi} for a rough discussion.

Below, we first present such a conjectural picture, motivated by the open question proposed by \cite{handel2022AFextremals} and an expected qualitative description of HL classes from our previous joint work with Shang \cite{hushangxiao2023hardlefArxiv}.

\subsection{Formulation of Shenfeld-van Handel's conjecture}
In this part, we formulate our conjectural picture on the numerical characterization of two-dimensional classes given by complete intersections of nef classes, refining the open proposal of Shenfeld and van Handel.

Before proceeding, we first introduce several notations and conventions that will be frequently used in the sequel.

Let $m$ be a positive integer and set $[m]= \{1,2,...,m\}$. Given vectors $v_1,...,v_m$ in a vector space $V$ and a nonempty subset $I\subset [m]$, we denote
$$v_I = \sum_{i\in I} v_i.$$
Given a real vector space $V$ and a subset $S \subset V$, we write $\Span_\mathbb{R} (S)$ for the linear subspace generated by $S$.

Let $X$ be a smooth projective variety of dimension $n$. We denote by $\prim(X)$ the space of prime divisors on $X$. The numerical dimension of a nef divisor or $(1,1)$-class $L$ is defined by $$\nd(L)=\max\{k \geq 0 | L^k \neq 0\}.$$
Let $m\leq n$ and let $\mathfrak{L}=(L_1,...,L_m)$ be a collection of nef classes on $X$. We use the same notation $\mathfrak{L}$ to denote the complete intersection class
$$\mathfrak{L}=L_1 \cdot...\cdot L_m. $$
For $J \subset [m]$, we write $\mathfrak{L}_{\setminus J}$ for the collection of nef classes $(L_k)_{k\in [m]\setminus J}$ or for the corresponding complete intersection class $$\mathfrak{L}_{\setminus J} =\prod_{k\in [m]\setminus J} L_k.$$

In this paper, we will mainly focus on the case when $m=n-2$.

The following terminology is taken from \cite{handel2022AFextremals}, where it is defined in the convexity setting.

\begin{defn}
A collection $\mathfrak{L}=(L_1,...,L_{m})$ of nef classes is called
\begin{itemize}
  \item subcritical, if $\nd(L_I) \geq |I|$ for any $I \subset [m]$;
  \item critical, if $\nd(L_I) \geq |I|+1$ for any $I \subset [m]$;
  \item supercritical, if $\nd(L_I) \geq |I|+2$ for any $I \subset [m]$.
\end{itemize}

\end{defn}

\begin{rmk}
By \cite{huxiaohardlef2022Arxiv} (see also Theorem \ref{nonvanishing}), the collection $\mathfrak{L}$ is subcritical if and only if $$\mathfrak{L}=L_1\cdot...\cdot L_{m}\neq 0.$$ 
Moreover, the class $\mathfrak{L}$ can be an HL class only if the collection is at least critical (see Proposition \ref{kernel nontriv}).
\end{rmk}

As we focus on the two dimensional case,
i.e., $m=n-2$, we denote $\ker \mathfrak{L}$ to be the subspace of divisor or $(1,1)$ classes $\eta$ such that $\mathfrak{L} \cdot \eta =0$.

\subsubsection{Expectation from positivity theory}
In \cite{hushangxiao2023hardlefArxiv}, by introducing a partial positivity notion -- ``$m$-lefness'' -- for algebraic maps via the defect of semismallness, we proved the following result (a special case of \cite{hushangxiao2023hardlefArxiv} when $p=q=1$). Let $X$ be a smooth projective variety of dimension $n$, and let $\mathfrak{L}=(L_1, ...,L_{n-2})$ be a collection of free line bundles on $X$. Assume that $L_I$ is $|I|+2$ lef for any $I\subset [n-2]$, then $\mathfrak{L}$ is an HL class on $H^{1,1}(X,\mathbb{R})$.

Recall that a free line bundle $L$ on $X$ is $m$-lef if for any irreducible subvariety $T$ on $X$,
\begin{equation*}
  2\dim T -\dim \Phi_L (T) -n \leq n-m,
\end{equation*}
where $\Phi_L: X \rightarrow \mathbb{P}(H^0 (X, L))$ is the Kodaira map associated with $L$.
This condition can be expressed numerically. In fact, the assumption that $L_I$ is $|I|+2$ lef can be formulated as follows:
for any irreducible subvariety $V$ such that $|I|+2-2\codim V >0$,
\begin{equation*}
  L_I ^{|I|+2-2\codim V} \cdot [V] \neq 0.
\end{equation*}
In particular, for any prime divisor $D$, $$L_I ^{|I|} \cdot [D] \neq 0.$$
By \cite{huxiaohardlef2022Arxiv}, the condition that $L_I ^{|I|} \cdot [D] \neq 0$ for every $I\subset [n-2]$ is equivalent to
\begin{equation*}
\mathfrak{L}\cdot [D]=L_1\cdot...\cdot L_{n-2}\cdot [D] \neq 0.
\end{equation*}

From a geometric viewpoint, at least when seeking for an HL class on $N^{1} (X)\subset H^{1,1}(X,\mathbb{R})$, it is natural to guess that this condition provides the main obstruction. On the other hand, by \cite{huxiaohardlef2022Arxiv} (see also \cite{Panov1987ONSP}), when $X$ is a compact complex torus of dimension $n$ and $\mathfrak{L}=(L_1, ...,L_{n-2})$ is a collection of nef classes on $X$, the class $\mathfrak{L}$ is an HL class on $H^{1,1}(X,\mathbb{R})$ if and only if $$\nd(L_I)\geq |I|+2$$ for every $I\subset [n-2]$, i.e., the collection $\mathfrak{L}$ is supercritical. These observations lead to the following qualitative expectation.

\begin{expec}\label{exp1}
Let $X$ be a smooth projective variety of dimension $n$, and let $\mathfrak{L}=(L_1, ...,L_{n-2})$ be a supercritical collection of free line bundles on $X$. Then $\mathfrak{L}$ is an HL class on $N^1 (X)$ if and only if
\begin{equation*}
\mathfrak{L}\cdot [D] =L_1\cdot...\cdot L_{n-2}\cdot [D] \neq 0.
\end{equation*}
for every prime divisor $D$ on $X$.
\end{expec}

\subsubsection{Expectation from convexity theory}
Another motivation comes from the recent breakthrough work of Shenfeld and van Handel \cite{handel2022AFextremals}. In that paper, the authors completely solved the long-standing open problem of characterizing of the extremals of the Alexandrov-Fenchel inequality for convex polytopes, dating back to Alexandrov's original paper \cite{alexandroff1938theorie} (see also \cite{alexandrovSelectedworksIconvexbodies} for an English translation). This characterization is essentially equivalent to describing ``HL classes of dimension two given by complete intersections of polytopes''.

Due to the deep connections of the Alexandrov-Fenchel inequality with many areas of mathematics, and the expectation that its algebraic formulation may extend to other settings, Shenfeld and van Handel proposed in \cite[Section 16.2]{handel2022AFextremals}
an algebraic analogue of their characterization of the extremals. 
We quote the following passage from Section 16.2 of their paper, which closely reflects our motivation:
\begin{quote}
``\emph{The rich algebraic theory surrounding the Alexandrov-Fenchel inequality raises the intriguing question whether our results might extend to a broader context.
This question arises, for example, if we aim to develop combinatorial applications
as in section 15 in situations that cannot be formulated in convex geometric terms...It may not be entirely obvious, however, how to even formulate algebraic analogues of the main results of this paper. The aim of this section is to sketch how our main results may be expressed in algebraic terms, which could (conjecturally) carry over to analogues of the Alexandrov-Fenchel inequality outside convexity. To the best of our knowledge, such problems are at present almost entirely open...}''

\end{quote}

Using toric geometry, Shenfeld and van Handel reformulated a special case of one of their main results \cite[Theorem 8.1]{handel2022AFextremals} as follows (see \cite[Corollary 16.2]{handel2022AFextremals}). Let $X_P$ be the toric variety corresponding to a simple polytope $P$ in $\mathbb{R}^n$ and let $\mathfrak{L}=(L_1,...,L_{n-2})$ be a collection of big and nef divisor classes on $X$, then any divisor class $\eta$ satisfying $\mathfrak{L}\cdot \eta=0$ can be written as $\eta=\eta_1 -\eta_2$ for some pseudo-effective divisor classes $\eta_1, \eta_2$ so that $\mathfrak{L}\cdot \eta_1 =\mathfrak{L}\cdot \eta_2 =0$. Regarding this result, Shenfeld and van Handel wrote:
\begin{quote}
``\emph{The point of the algebraic formulation is, however, that
the same algebraic structures carry over to other mathematical problems \cite{huhICM2018}. The
statement of Corollary 16.2 (for example) therefore gives rise to natural conjectures on what analogues of the results of this paper might look like in other contexts...One might therefore ask whether a result such as Corollary 16.2 carries over to this setting, at least in sufficiently nice situations. To the best of our knowledge this question is entirely open, except for toric varieties which admit a precise correspondence with convex geometry \cite{fultonToricBook, ewaldConvexity} (for which such a conclusion follows from the results of this paper)...}''
\end{quote}

A natural generalization of \cite[Corollary 16.2]{handel2022AFextremals} to arbitrary algebraic varieties leads to the following expectation.

\begin{expec}\label{exp2}
Let $X$ be a smooth projective variety of dimension $n$, and let $\mathfrak{L}=(L_1, ...,L_{n-2})$ be a collection of big nef classes on $X$. Then any element $\eta\in\ker \mathfrak{L}$ can be written as $\eta=\eta_1 -\eta_2$ for some pseudo-effective classes $\eta_1, \eta_2$ so that $\mathfrak{L}\cdot \eta_1 =\mathfrak{L}\cdot \eta_2 =0$.
\end{expec}

The expectation provides a quantitative description of HL classes.
Note that on a toric variety, every pseudo-effective class is effective. 
Combining Expectations \ref{exp1}, \ref{exp2} and the full statement of \cite[Theorem 8.1]{handel2022AFextremals}, we propose the following conjectural picture in the supercritical case, as suggested by Shenfeld-van Handel:

\begin{conj}[Analytic supercritical collections]\label{main conj}
Let $X$ be a compact K\"ahler manifold of dimension $n$, and let $\mathfrak{L}=(L_1,...,L_{n-2})$ be a supercritical collection of nef classes on $X$. Define the vector spaces $V_{\mathfrak{L},\eff}$ as follows:
\begin{align*}
V_{\mathfrak{L}, \eff} = \Span_{\mathbb{R}} \{[D]: D\in \prim(X),\ \mathfrak{L} \cdot [D]=0\}.
\end{align*}
Then $$\ker (\mathfrak{L}:H^{1,1}(X,\mathbb{R})\rightarrow H^{n-1,n-1}(X,\mathbb{R})) = V_{\mathfrak{L},\eff}.$$
In particular, $\mathfrak{L}$ is an HL class on $H^{1,1}(X,\mathbb{R})$ if and only if $\mathfrak{L} \cdot [D] \neq 0$ for every prime divisor $D$ on $X$.
\end{conj}

\begin{conj}[Algebraic supercritical collections]\label{main conj alg}
  Let $X$ be a smooth projective variety over an algebraically closed field $k$ of dimension $n$, and let $\mathfrak{L}=(L_1,...,L_{n-2})$ be a supercritical collection of nef divisors on $X$. 
  Then $$\ker (\mathfrak{L}:N^1(X)\rightarrow N^{n-1}(X)) = V_{\mathfrak{L},\eff}.$$
  In particular, $\mathfrak{L}$ is an HL class on $N^1(X)$ if and only if $\mathfrak{L} \cdot [D] \neq 0$ for every prime divisor on $X$.
\end{conj}

As our main progress concerns the supercritical case, we provide in this introduction only a brief discussion of the conjectural picture for arbitrary collections. The general case is considerably more subtle and requires new ideas; further details are referred to Section \ref{sec arbi}. Our formulation is inspired by \cite[Theorem 2.13]{handel2022AFextremals}. We begin with the following definition.

\begin{defn}
Let $\mathfrak{L}=(L_1,...,L_{n-2})$ be a collection of nef classes and fix an ample class $A$ on $X$. A pair $(\alpha, \beta)$ of nef classes is called a $\mathfrak{L}$-degenerate pair if
\begin{align*}
\mathfrak{L}\cdot {\alpha} \cdot {\beta} =0\text{ and } \mathfrak{L}\cdot A \cdot {\alpha}=\mathfrak{L} \cdot A \cdot {\beta }.
\end{align*}
\end{defn}

\begin{rmk}
By \cite{huxiaohardlef2022Arxiv} (see also Theorem \ref{prop geom}), the condition on the pair can be replaced by
\begin{align*}
\mathfrak{L}\cdot {\alpha} \cdot {\beta} =0,\ \mathfrak{L}\cdot  {(\alpha-\beta)}=0.
\end{align*}
Hence the notion of degenerate pairs is independent of the choice of the ample class $A$.
\end{rmk}

\begin{conj}[Arbitrary analytic collections]\label{main conj1}
Let $X$ be a compact K\"ahler manifold of dimension $n$, and let $\mathfrak{L}=(L_1,...,L_{n-2})$ be a collection of nef classes on $X$. Define the vector space $V_{\mathfrak{L},\eff}$ as in the supercritical case, and set
\begin{align*}
V_{\mathfrak{L},\deg}={\Span}_{\mathbb{R}}\{\mu_*(\widetilde{\alpha}-\widetilde{\beta})| (\widetilde{\alpha}, \widetilde{\beta})\ \text{is a $\mu^*\mathfrak{L}$-degenerate pair on $\widetilde{X}$}\},
\end{align*}
where the span is taken over all K\"ahler modifications $\mu : \widetilde{X}\rightarrow X$ and all $\mu^*\mathfrak{L}$-degenerate pairs $(\widetilde{\alpha}, \widetilde{\beta})$ on $\widetilde{X}$.
Then $$\ker (\mathfrak{L}:H^{1,1}(X,\mathbb{R})\rightarrow H^{n-1,n-1}(X,\mathbb{R})) = V_{\mathfrak{L},\eff}+V_{\mathfrak{L},\deg}.$$
In particular, $\mathfrak{L}$ is an HL class if and only if $V_{\mathfrak{L},\deg}=\{0\}$ and $\mathfrak{L} \cdot [D] \neq 0$ for every prime divisor $D\subset X$.
\end{conj}

\begin{rmk}
  We also conjecture that $\ker (\mathfrak{L}:N^1(X)\rightarrow N^{n-1}(X))$ has    
  the same structure as predicted in Conjecture \ref{main conj1} for any smooth projective variety over an algebraically closed field. 
\end{rmk}

\begin{rmk}
When the collection $\mathfrak{L}$ is not subcritical, by \cite{huxiaohardlef2022Arxiv} the complete intersection class $$\mathfrak{L}=L_1\cdot...\cdot L_{n-2}= 0,$$
so that $\ker \mathfrak{L}$ coincides with the whole space and Conjecture \ref{main conj1} holds trivially in this case. Therefore, the conjecture is nontrivial only when the collection is subcritical.
It is straightforward to check that $V_{\mathfrak{L},\deg} \subset \ker \mathfrak{L}$, and moreover, when $\mathfrak{L}$ is supercritical, we have $V_{\mathfrak{L},\deg}=\{0\}$ (see Proposition \ref{vdeg}).
\end{rmk}

\begin{rmk}
Regarding the space $V_{\mathfrak{L},\eff}$, we will show that whenever the collection $\mathfrak{L}$ is critical, the number of prime divisors $D$ satisfying $\mathfrak{L}\cdot [D]=0$ is finite. Indeed, the classes $[D]$ must be linearly independent in the real N\'{e}ron-Severi space, so their number is at most the Picard number of $X$. See Section \ref{sec localaf} for details.
\end{rmk}

\begin{rmk}
Although Conjectures \ref{main conj} and \ref{main conj1} are formulated for smooth projective varieties and K\"ahler manifolds, the statements naturally extend to other mathematical objects where analogues of intersection numbers, nef classes, effective classes and birational modifications can be defined.
For instance, this framework applies to fan varieties.
\end{rmk}

\subsection{The main results}

In the sequel, we make the following conventions: otherwise stated, the projective variety is defined over an algebraically closed field of characteristic 0; for a compact K\"ahler manifold $X$, we denote by $\ker \mathfrak{L}$ the kernel of $\mathfrak{L}$ acting on $H^{1,1}(X,\mathbb{R})$; for a projective variety, however, the same symbol $\ker \mathfrak{L}$ denotes the kernel of $\mathfrak{L}$ acting on $N^1(X)$.

We focus primarily on the supercritical setting and obtain the following result:

\begin{thrmx}\label{ker supercritical intro}
Let $X$ be a smooth projective variety of dimension $n$, and let $\mathfrak{L}=(L_1,...,L_{n-2})$ be a collection of nef classes on $X$. Assume that the collection $\mathfrak{L}$ is supercritical and satisfies 
$$\nd(L_k)\geq k+2,\ \forall\ k\in[n-2],$$
then Conjecture \ref{main conj alg} holds, i.e., $\ker \mathfrak{L} = V_{\mathfrak{L},\eff}.$

As a consequence, given an arbitrary full (increasing) flag $\mathcal{F}$ of subsets of $[n-2]$, that is, a sequence
\begin{equation*}
  \mathcal{F}: I_1 \subset I_2 \subset...\subset I_{n-2} =[n-2],
\end{equation*}
with $|I_k| =k$, and an arbitrary supercritical collection $\mathfrak{L}$, one can construct a new collection of nef classes
\begin{equation*}
  \mathfrak{L}_\mathcal{F}=(L_{I_1},...,L_{I_{n-2}}),
\end{equation*}
for which Conjecture \ref{main conj alg} holds:
$$\ker\mathfrak{L}_\mathcal{F}=
V_{\mathfrak{L}_\mathcal{F},\eff}.$$ Moreover, the following equality always holds:
\begin{equation*}
  \sum_{\mathcal{F}} \ker \mathfrak{L}_\mathcal{F}  =V_{\mathfrak{L},\eff},
\end{equation*}
where the sum is taken over all full flags $\mathcal{F}$ of $[n-2]$.

\end{thrmx}

Regarding the requirement $\nd(L_k)\geq k+2 $ for each $k\in[n-2]$, this assumption is needed in our argument because the proof of the characterization of $\ker \mathfrak{L}$ proceeds inductively  using the algebraic analogue of Shenfeld-van Handel's local Alexandrov-Fenchel inequality (see Theorem \ref{local af intr} below). Roughly speaking, the assumption guarantees bigness at each step of the induction. For a more accessible illustration of the idea, see Theorem \ref{ker bignef}, which treats the case where every class in the collection is big and nef.

\begin{rmk}
By Theorem \ref{ker supercritical intro}, the description of $\ker \mathfrak{L}$ in Conjecture \ref{main conj} for a general supercritical collection is equivalent to
\begin{equation*}
  \sum_{\mathcal{F}} \ker \mathfrak{L}_\mathcal{F} =\ker \mathfrak{L}.
\end{equation*}
\end{rmk}

\begin{exmple}
As a typical example, if every $L_k$ in the collection $\mathfrak{L}$ is big and nef and the augmented base locus (equivalently, the null locus) of each $L_k$ has codimension at least 2, then $\mathfrak{L} \cdot [D]\neq 0$ for any prime divisor $D\subset X$. Consequently, the complete intersection $$\mathfrak{L}=L_1 \cdot...\cdot L_{n-2}$$ is an HL class. The same statement holds on a compact K\"ahler manifold.
To our best knowledge, this provides the first series of two-dimensional HL classes arising as complete intersections of nef classes for which neither freeness nor semi-positivity is required.
\end{exmple}

As an application, we apply Theorem \ref{ker supercritical intro} to characterize the extremals of log-concave sequence defined by the intersection numbers of big nef classes (see Section \ref{sec hodge extre} for a more general version).

\begin{corx}\label{corx extrem logc}
Let $X$ be a smooth projective variety of dimension $n$, and let $A, B$ be big and nef classes on $X$.
Define $$a_l = A^l \cdot B^{n-l},\ 0\leq l\leq n.$$
Fix $k\in [n-1]$. Then the equality
$a_k^2 = a_{k-1} a_{k+1}$
holds if and only if
\begin{equation*}
  A-cB =\sum_{i} c_i [D_i]
\end{equation*}
for some prime divisors $D_i$ satisfying $$A^{k-1} \cdot B^{n-k-1}\cdot [D_i]=0.$$
Here, $c>0$ is the constant satisfying $$A^{k} \cdot B^{n-k-1}=cA^{k-1} \cdot B^{n-k}$$ whose existence is guaranteed by Theorem \ref{prop geom}.
\end{corx}

Note that the sequence $\{a_l\}_l$ is always log-concave, i.e, $$a_l ^2 \geq a_{l-1} a_{l+1}, \forall 1\leq l\leq n-1.$$ When equality holds for every $l$, this is Teissier's proportionality problem (see \cite{BFJ09, fuxiaoTeissProp}). The extremals of $a_k^2 = a_{k-1} a_{k+1}$ for a fixed $k$ had remained open. Corollary \ref{corx extrem logc} resolves this question for $A, B$ big and nef.

Our approach towards the conjectures is inspired by \cite{handel2022AFextremals}, where Shenfeld and van Handel proved a local Alexandrov-Fenchel inequality that enables an inductive argument. In our setting, we prove the following local Hodge index inequality, which serves as the algebraic analogue of the local Alexandrov-Fenchel inequality.

\begin{thrmx}\label{local af intr}
Let $X$ be a smooth projective variety of dimension $n$, and
let $\mathfrak{L}=(L_1,...,L_{n-2})$ be a critical collection of nef classes on $X$. Fix some $r\in [n-2]$.
Then for any $\alpha \in \ker \mathfrak{L}$, there exists a class $\beta$ such that
\begin{enumerate}
  \item $\beta - \alpha \in V_{\mathfrak{L},\eff}$;
  \item $-\beta^2 \cdot \mathfrak{L}_{\setminus r} \in \Mov_1 (X)$,
\end{enumerate}
where $\Mov_1 (X)$ denotes the closed convex cone generated by movable curve classes.
\end{thrmx}

More generally, we establish a broader statement  -- the local Hodge index inequality for Lorentzian polynomials -- which contains Theorem \ref{local af intr} as a special case (see Theorem \ref{local af}). The proof abstracts the scheme of Shenfeld-van Handel. However, in the algebro-geometric and analytic-geometric versions and their applications, we need to construct a generating set of the pseudo-effective cone with certain special properties (see Section \ref{sec generating set}). Achieving this relies crucially on the positivity theory of divisor and $(1,1)$-classes developed in the last decades, including divisorial Zariski decompositions \cite{Bou04,Nak04}, intersection theory of pseudo-effective classes \cite{BFJ09,BDPP13, lehmannXiaoPosiCurve}, geometric properties of the null locus and the augmented base locus \cite{BFJ09, nakamayeNullLocus,elmnp09restrictedvolume, tosattiNullLocus, nystromDualityMorse}.

In this framework, the local Hodge index inequality applies not only to divisor classes on projective manifolds, but also to transcendental classes on complex projective manifolds and to Lorentzian fans in full generality. Moreover, it extends to compact K\"ahler manifolds when every element in $\mathfrak{L}=(L_1,...,L_{n-2})$ is big and nef (see Section \ref{sec localaf} for more details).
This result allows us to extend Theorem \ref{ker supercritical intro} to transcendental classes on complex projective manifolds in full generality, as well as to compact K\"ahler manifolds or to Bergman fans of a matroids and polymatroids under appropriate positivity assumptions.

\subsection*{Organization}
The paper is organized as follows. In Section \ref{sec prel}, we give a brief introduction to positivity in geometry, Lorentzian polynomials and Lorentzian fans. In Section \ref{sec posicri}, we recall some of our previous results and extend the notion of numerical dimension and the criterion for nonvanishing intersections to Lorentzian polynomials. In Section \ref{sec collec}, we present several basic properties of collections of nef classes. Section \ref{sec localaf} is devoted to the proof of local Hodge index inequality and to the construction of generating sets in various contexts. In Section \ref{sec char}, we give the numerical characterization of hard Lefschetz classes. Finally, in Section \ref{sec arbi} we discuss in more details on the conjectural picture for an arbitrary collection and provide several examples.

\subsection*{Acknowledgements}
This work is supported by the National Key Research and Development Program of China (No. 2021YFA1002300) and National Natural Science Foundation of China (No. 11901336).
We would like to thank Junyan Cao, Ya Deng and Jie Liu for helpful discussions. In particular, we would like to thank Junyan Cao for the help of Lemma \ref{union locus}. We also thank the referee for helpful comments.
The impact of Shenfeld-van Handel's breakthrough work on this article cannot be overestimated, and we are very grateful to Ramon van Handel for his kind messages and numerous helpful comments and questions.
We would like to point it out that the algebraic analogue of the local Alexandrov-Fenchel inequality was also realized by Karim Adiprasito and Ramon van Handel.

\section{Preliminaries}\label{sec prel}

In this section, we introduce some notions and present several preliminary results which will be applied in the sequel.

\subsection{Positivity}\label{sec pos}
We first briefly introduce some positivity notions in algebraic geometry and analytic geometry. For more details, we refer the reader to \cite{dem_analyticAG, Dem_AGbook} for the analytic setting and \cite{lazarsfeldPosI} for the algebraic setting.

\subsubsection{Algebraic positivity}
\begin{defn}
Let $X$ be a smooth projective variety of dimension $n$.
\begin{itemize}
\item $N^1(X)$ is the real vector space of numerical classes of divisors.
\item $N_1(X)$ is the real vector space of numerical classes of curves.
\item $\prim(X)$ is the set of prime divisors on $X$.
\item $\Eff^1(X)$ is the cone of pseudo-effective divisor classes, that is, an element $L\in \Eff^1(X)$ can be written as $$L=\lim_k [F_k],$$ where $F_k$ is a non-negative combination of elements in $\prim(X)$. It is clear that $\Eff^1(X)$ is a closed convex cone in $N^1(X)$. The interior of $\Eff^1(X)$ is called the big cone of  divisor classes, and an interior point is called a big divisor class.
\item $\Nef^1(X)$ is the cone of nef divisor classes, that is, $L\in \Nef^1(X)$ if and only if $L\cdot C\geq 0$ for any irreducible curve $C$ on $X$.
The interior of $\Nef^1(X)$ is the ample cone $\Amp^1 (X)$, and an interior point is called an ample class.
\item $\Mov^1(X)$ is the cone of movable divisor classes, that is, $L\in \Mov^1(X)$ if and only if $$L = \lim_{k\rightarrow \infty} (\pi_k)_* A_k,$$ where $\pi_k: X_k \rightarrow X$ is a birational morphism and $A_k \in \Amp^1 (X_k)$ is a real ample divisor class on $X_k$.
\item $\Mov_1(X)$ is the cone of movable curve classes, that is, $\Mov_1(X)$ is the closure of the set of convex combinations of curve classes of the form $$C = \lim_{k\rightarrow \infty} (\pi_k)_* (A_1 \cdot...\cdot A_{n-1}),$$
 where $\pi_k: X_k \rightarrow X$ is a birational morphism and the $A_i \in \Amp^1 (X_k)$ are real ample divisor classes on $X_k$.
\end{itemize}

In the above algebraic setting, all the cones are either in $N^1(X)$ or $N_1(X)$. On a projective surface, we have that: $\Nef^1(X)=\Mov^1(X)=\Mov_1(X)$.

\end{defn}

\begin{rmk}
By the fundamental result in \cite{BDPP13}, on a smooth projective variety $X$, the dual cone of $\Eff^1(X)$ is exactly given the cone of movable curve classes:
\begin{equation*}
  \Eff^1(X) ^* = \Mov_1(X).
\end{equation*}
\end{rmk}

\subsubsection{Analytic positivity}
The algebraic positivity notions have their analytic counterparts when one studies the analytic geometry of a smooth complex projective variety or a possibly non-projective compact K\"ahler manifold.

\begin{defn}
Let $X$ be a compact K\"ahler manifold of dimension $n$.
\begin{itemize}
\item $H^{1,1}(X, \mathbb{R})$ is the real vector space of $(1,1)$ classes.
\item $H^{n-1,n-1}(X, \mathbb{R})$ is the real vector space of $(n-1,n-1)$ classes.
\item $\Eff^1(X)$ is the cone of pseudo-effective $(1,1)$ classes, that is, $L \in \Eff^1(X)$ if and only if it can be represented by a positive $(1,1)$ current. An interior point of this cone is called a big $(1,1)$ class.
\item $\Nef^1(X)$ is the cone of nef $(1,1)$ classes. Its interior is the K\"ahler cone, consisting of K\"ahler classes.
\item $\Mov^1(X)$ is the cone of movable $(1,1)$ classes: $L\in \Mov^1(X)$ if and only if $$L = \lim_{k\rightarrow \infty} (\pi_k)_* \omega_k,$$ where $\pi_k: X_k \rightarrow X$ is a K\"ahler modification and $\omega_k$ is a K\"ahler class on $X_k$.
\item $\Mov_1(X)$ is the cone of movable $(n-1,n-1)$ classes, that is, $\Mov_1(X)$ is the closure of the set of convex combinations of classes of the form $$C = \lim_{k\rightarrow \infty} (\pi_k)_* (\omega_1 \cdot...\cdot \omega_{n-1}),$$
 where $\pi_k: X_k \rightarrow X$ is a K\"ahler modification and the $\omega_i$ are K\"ahler classes on $X_k$.
\end{itemize}

\end{defn}

In the analytic setting, all the cones are either in $H^{1,1}(X, \mathbb{R})$ or $H^{n-1,n-1}(X, \mathbb{R})$.

\begin{rmk}
The K\"ahler analog of \cite{BDPP13} regarding the duality of cones
\begin{equation*}
  \Eff^1(X) ^* = \Mov_1(X)
\end{equation*} 
still remains open when $X$ is an arbitrary compact K\"ahler manifold.
Nevertheless, by \cite{nystromDualityMorse}, it
holds in the analytic setting when $X$ is a smooth complex projective variety.
\end{rmk}

\emph{\textbf{Convention}}: Because of their similarities, both in the algebraic and analytic settings we have already used the same symbols to denotes the corresponding positive cones. In the sequel, we make the convention: unless otherwise specified, when a result holds in both settings, we shall use the same notations for these cones.

\subsubsection{Zariski decomposition}\label{sec zardeco}
By Nakayama \cite{Nak04} for the algebraic setting and Boucksom \cite{Bou04} for a general compact complex manifold, the cones $\Eff^1(X), \Mov^1(X), \Nef^1 (X)$ are closely related via divisorial Zariski decomposition. More precisely, when $X$ is a projective manifold or a compact K\"ahler manifold, every $L\in \Eff^1(X)$ can be decomposed as follows:
\begin{equation*}
  L=P(L)+[N(L)],
\end{equation*}
where $P(L) \in \Mov^1(X)$ and $N(L)$ is an effective divisor determined by the class $L$. Furthermore, $N(L)$ is supported by at most $\rho(X)$ prime divisors, where $\rho(X)$ is the Picard number of $X$.

\subsubsection{Augmented base locus and null locus}
We recall some basic properties on the augmented base locus and the null locus of a big nef class.

\begin{defn}
Let $X$ be a smooth projective variety.
Given a big real divisor class $L$ on $X$, its augmented base locus (or non-ample locus) is defined by
\begin{equation*}
  \mathbb{B}_+ (L) = \bigcap_{L=A+E} \supp E,
\end{equation*}
where the intersection is taken over all decompositions $L=A+E$ such that $A$ is a real ample divisor and $E$ is a real effective divisor.

In the analytic setting, let $X$ be a compact K\"ahler manifold and $L$ a big $(1,1)$ class, then its augmented base locus (or non-K\"ahler locus) is defined by
\begin{equation*}
  \mathbb{B}_+ (L) = \bigcap_{T \in L} \sing T,
\end{equation*}
where the intersection is taken over all K\"ahler currents $T\in L$ with analytic singularities and $\sing T$ is the set of singularities of $T$.

\end{defn}

It is clear that $\mathbb{B}_+ (L)$ is a proper subvariety for $L$ big.
By \cite{Bou04}, the non-K\"ahler locus and non-ample locus coincide when $L$ is a divisor class.

\begin{defn}
Let $L$ be a nef class, the null locus of $L$ is defined by
\begin{equation*}
  \Null(L)=\bigcup_{\ L^{\dim V} \cdot [V]=0}  V,
\end{equation*}
where the union is taken over all subvarieties $V \subset X$ such that $L^{\dim V} \cdot [V]=0$.
\end{defn}

The augmented base locus is closely related to the null locus.

\begin{lem}
[see \cite{nakamayeNullLocus, elmnp09restrictedvolume, tosattiNullLocus}]\label{augmented null locus}
Let $X$ be a smooth projective variety (or in the analytic setting, a compact K\"ahler manifold) of dimension $n$.
For a nef class $L$ on $X$, we have that $$\mathbb{B}_+ (L)=\Null(L).$$
In particular, $\Null(L)$ is proper when $L$ is big.
\end{lem}

\subsubsection{Extension to other fields}
\label{fields not c}
While we only work on varieties over an arbitrary algebraically closed field of characteristic 0. The results mentioned in this section extend to smooth projective varieties over an arbitrary algebraically closed field, we need the following necessary generalizations, which are verified in the references:
\begin{itemize}
\item The paper \cite{birkarNullLocusAnyfield} proves that the augmented base locus and the null locus coincide over an arbitrary field.
\item The paper \cite{fulgerlehmZariskiCycle} describes how to extend \cite{BDPP13, BFJ09} to an arbitrary algebraically closed field.
\end{itemize}

\subsection{Lorentzian polynomials}
The theory of Lorentzian polynomials was introduced and systematically developed in \cite{branhuhlorentz} and independently (with part overlap) in \cite{logconcavepoly1, logconcavepoly2, logconcavepoly3}. The class of Lorentzian polynomials contains all homogeneous stable polynomials, and is intimately connected to matroid theory, negative dependence properties, Potts model partition functions and log-concave polynomials. Since the volume polynomials of nef divisors on a projective variety and the volume polynomials of convex bodies are Lorentzian, they also reveal important information on projective varieties and convex bodies (see \cite{huh2022combinatorics} for a nice exposition).
The class of Lorentzian polynomials can be considered as an analog of the Hodge-Riemann relation of degree one, it has many remarkable applications in combinatorics and geometry.

In the purely ``polynomial proof'' of the Heron-Rota-Welsh conjecture that does not rely on Hodge theory, Br\"{a}nd\'{e}n-Leake \cite{branden2023lorentziancone} introduced a more general notion of Lorentzian polynomials on cones.

\begin{defn}
Let $\mathfrak{C}$ be an open convex cone in $\mathbb{R}^s$. A homogeneous polynomial  $f$ (with $\deg f=n$) on $\mathbb{R}^s$ is called $\mathfrak{C}$-Lorentzian if for all $v_1, ...,v_n \in \mathfrak{C}$,
\begin{itemize}
  \item $D_{v_1}... D_{v_n} f >0$, and
  \item the symmetric bilinear form $$(\xi, \eta) \mapsto D_\xi D_\eta D_{v_3}... D_{v_n} f$$ has exactly one positive eigenvalue.
\end{itemize}
Here, $D_v$ is the directional derivative along $v$.
\end{defn}

When $\mathfrak{C}=\mathbb{R}_{>0} ^s$, these are the Lorentzian polynomials studied in \cite{branhuhlorentz}.
The above definition is equivalent to that for all positive integers $m$ and for all $v_1,...,v_m \in \mathfrak{C}$, the polynomial
\begin{equation*}
  (y_1,...,y_m)\mapsto f(y_1 v_1 +...+ y_m v_m)
\end{equation*}
is Lorentzian in the sense of \cite{branhuhlorentz} and has only positive coefficients.

\subsection{Lorentzian fans}\label{sec lorfan}
From the perspective of fans, Dustin Ross \cite{ross2023lorentzianfan} introduced a very interesting geometric counterpart for Lorentzian polynomials -- Lorentzian fans (see also \cite{branden2023lorentziancone} from the perspective of hereditary polynomials).
In the following, we introduce several notions and results on Ross's construction.

Let $\Sigma$ be a (finite) collection of strongly convex polyhedral cones in a vector space $V$. Given two cones $\tau, \sigma$, $\tau \preceq \sigma$ means that $\tau$ is a face of $\sigma$.
The collection $\Sigma$ is a fan if it satisfies the following properties:
\begin{itemize}
  \item if $\sigma \in \Sigma$ and $\tau \preceq \sigma$, then $\tau \in \Sigma$;
  \item if $\sigma_1, \sigma_2 \in \Sigma$, then $\sigma_1 \cap \sigma_2 \preceq \sigma_1$ and $\sigma_1 \cap \sigma_2 \preceq \sigma_2$.
\end{itemize}

The collection of $k$-dimensional cones of $\Sigma$ is denoted by $\Sigma(k)$. In particular, $\Sigma(1)$ is the set of rays. The fan $\Sigma$ is called pure if every inclusion-maximal cone of $\Sigma$ has the same dimension, and it is called a $n$-fan if it is pure of dimension $n$. The support of $\Sigma$, $|\Sigma|$, is the union of all cones $\sigma \in \Sigma$.

The fan $\Sigma$ is called
simplicial if all cones are simplicial, i.e., for any $\sigma\in \Sigma$, $\dim \sigma = |\sigma (1)|$, where $\sigma (1)$ is the set of rays in $\sigma$. In the sequel, we always assume that $\Sigma$ is simplicial.

For any ray $\rho \in \Sigma (1)$, we fix a non-zero marketing vector $u_{\rho}$ on $\rho$ (for example, on a rational fan one can take $u_{\rho}$ to be the unique primitive integral vector on the ray $\rho$).

Given a simplicial fan $\Sigma$ in $V$, the Chow ring of $\Sigma$ is defined by
\begin{equation*}
  A^{\bullet} (\Sigma) =\frac{\mathbb{R}[x_\rho| \rho \in \Sigma(1)]}{\mathcal{I}+\mathcal{J}},
\end{equation*}
where
\begin{equation*}
  \mathcal{I}=\langle x_{\rho_1}... x_{\rho_k}| \cone(\rho_1,...,\rho_k) \notin \Sigma\rangle,\ \mathcal{J}=\langle \sum_{\rho\in \Sigma (1)} (v, u_{\rho}) x_{\rho}| v \in V^{*}\rangle,
\end{equation*}
where $V^*$ is the dual space of $V$, and the bracket $(-, -)$ denotes the pairing between $V^*$ and $V$.
Note that both $\mathcal{I}$ and $\mathcal{J}$ are homogeneous, the Chow ring $A^{\bullet} (\Sigma)$ is a graded ring. Denote by $A^k (\Sigma)$ the subgroup of homogeneous elements of degree $k$. The Chow ring $A^{\bullet} (\Sigma)$ is generated by elements of degree one.
Denote the generators of $A^{\bullet} (\Sigma)$ by
\begin{equation*}
  X_\rho =[x_\rho]\in A^{1} (\Sigma),
\end{equation*}
and for any $\sigma \in \Sigma(k)$, define the cone monomial
\begin{equation*}
  X_\sigma=\prod_{\rho\in \sigma(1)} X_\rho \in A^{k} (\Sigma).
\end{equation*}
By \cite{huhHRR}, if $\Sigma$ is a simplicial $n$-fan, then $A^k(\Sigma)$ is spanned by $X_\sigma$ with $\sigma\in \Sigma(k)$ and $A^k (\Sigma)=0$ when $k>n$.

Given $\tau\in \Sigma$, denote $V_\tau$ the subspace generated by the cone $\tau$. Recall that a weight function $\omega: \Sigma(k) \rightarrow \mathbb{R}$ satisfying the weighted balanced condition:
\begin{equation*}
  \sum_{\sigma\in \Sigma(k), \tau \preceq \sigma} \omega(\sigma) u_{\sigma\setminus \tau} \in V_{\tau},\ \forall \tau \in \Sigma(k-1),
\end{equation*}
is called a Minkowski $k$-weight. Here, ${\sigma\setminus \tau}$ is the unique ray in $\sigma(1)\setminus \tau(1)$.

If $\Sigma$ is a simplicial $n$-fan endowed with a positive Minkowski $n$-weight $\omega$, then the pair $(\Sigma, \omega)$ is called a tropical fan. It is called a balanced fan if the weight $\omega \equiv 1$ is a Minkowski $n$-weight.
By \cite{huhHRR}, when $(\Sigma, \omega)$ is tropical, there is a well-defined degree map
\begin{equation*}
  \deg_\Sigma: A^n(\Sigma) \rightarrow \mathbb{R}
\end{equation*}
such that $\deg(X_\sigma) = \omega(\sigma)$ for every $\sigma\in \Sigma(n)$.

Given a cone $\tau\in \Sigma$, the neighborhood of $\tau$ in $\Sigma$ is defined by
\begin{equation*}
  N_\tau \Sigma =\{\pi| \pi\preceq \sigma\ \text{for some}\ \sigma\in \Sigma\ \text{with}\ \tau\preceq \sigma\}.
\end{equation*}
Let $V_\tau$ be the subspace generated by $\tau$. The star fan of $\Sigma$ at $\tau$ is the fan in the quotient space $V^\tau = V/ V_\tau$ defined by the quotients of cones in $ N_\tau \Sigma$:
\begin{equation*}
  \Sigma^\tau = \{\overline{\pi} \subset V^\tau | \pi\in  N_\tau \Sigma \}.
\end{equation*}
It can be checked that $\Sigma^\tau$ is a simplicial $n^\tau$-fan ($n^\tau = n- \dim V_\tau$) when $\Sigma$ is a simplicial $n$-fan, and that the marketing vectors $u_\rho$ induce a marketing on the quotient fan. Moreover, by \cite{ross2023lorentzianfan}, a Minkowski $n$-weight $\omega$ determines a Minkowski $n^\tau$-weight $\omega^\tau$ on $\Sigma^\tau$. The determination is unique up to some positive scaling. In particular, $(\Sigma^\tau, \omega^\tau)$ is tropical whenever $(\Sigma, \omega)$ is tropical.

The vector space of divisor classes on $\Sigma$ is defined by
\begin{equation*}
  D(\Sigma)=\frac{PL(\Sigma)}{L(\Sigma)},
\end{equation*}
where $PL(\Sigma)$ is the space of piecewise linear functions on $\Sigma$, that is, $\phi: |\Sigma|\rightarrow \mathbb{R}$ is in $PL(\Sigma)$ if for any $\sigma\in \Sigma$, there is a linear function $\phi_\sigma \in V^*$ such that $\phi_{|\sigma}=\phi_\sigma$, and $L(\sigma)\subset PL(\Sigma)$ is the subspace given by the restrictions of linear functions $\phi \in V^*$. Let $D_\rho$ be the class of the piecewise linear function taking value $1$ at $u_{\rho}$ and value $0$ at $u_\eta$ for $\eta\neq \rho$, then $\{D_\rho| \rho\in \Sigma(1)\}$ spans $D(\Sigma)$ and
\begin{equation*}
  D(\Sigma)\rightarrow A^1(\Sigma),\ D_\rho \mapsto X_\rho,
\end{equation*}
gives an isomorphism of vector spaces. The volume function on $D(\Sigma)$ is defined by
\begin{equation*}
\vol_\Sigma (D) =\deg_\Sigma (D^n),\ D\in D(\Sigma),
\end{equation*}
by identifying $D$ as an element of $A^1 (\Sigma)$.

\subsubsection{Positive cones on a fan}
The pseudo-effective cone of divisors on $\Sigma$ is the polyhedral cone generated by $\{D_\rho| \rho\in \Sigma(1)\}$:
\begin{equation*}
  \Eff^1(\Sigma)=\cone\{D_\rho| \rho\in \Sigma(1)\},
\end{equation*}
that is, $\alpha\in \Eff^1(\Sigma)$ if and only if there are $c_\rho \geq 0$ such that $\alpha = \sum_{\rho \in \Sigma(1)} c_\rho D_\rho$.

Following \cite{ross2023lorentzianfan} and \cite{huhHRR},
\begin{equation*}
  \Nef^1(\Sigma)=\{\alpha\in D(\Sigma)|\ \forall \tau \in \Sigma,\ \exists\
  \text{a representative}\ \phi\ \text{such that}\ \phi_{|\tau}=0,\ \phi_{|N_\tau \Sigma}\geq 0\},
\end{equation*}
\begin{equation*}
  \Amp^1(\Sigma)=\{\alpha\in D(\Sigma)|\ \forall \tau \in \Sigma,\ \exists\
  \text{a representative}\ \phi\ \text{such that}\ \phi_{|\tau}=0,\ \phi_{|N_\tau \Sigma \backslash \tau}> 0\}.
\end{equation*}

Then the interior of $\Nef^1(\Sigma)$ is given by $\Amp^1(\Sigma)$.

The fan $\Sigma$ is called quasi-projective if $\Amp^1(\Sigma)\neq \emptyset$.

\begin{defn}
A tropical $n$-fan $(\Sigma, \omega)$ is called Lorentzian if $\Sigma$ is quasi-projective and for any $\tau\in \Sigma$ and all $D_3,...,D_{n^\tau}\in \Amp^1(\Sigma^\tau)$, the quadratic form
\begin{align*}
  &D(\Sigma^\tau)\times D(\Sigma^\tau)\rightarrow \mathbb{R},\\
  &(D_1, D_2)\mapsto \deg_{\Sigma^\tau} (D_1\cdot D_2 \cdot D_3\cdot...\cdot D_{n^\tau})
\end{align*}
has exactly one positive eigenvalue.
\end{defn}

\begin{prop}\label{lorfan}
Let $\Sigma$ be a quasi-projective tropical $n$-fan, then it is Lorentzian if and only if the volume polynomial $\vol_{\Sigma^\tau}$ is $\Amp^1(\Sigma^\tau)$-Lorentzian for any $\tau\in \Sigma$.
\end{prop}


\subsubsection{Bergman fan of a matroid}
An important class of Lorentzian fans is the class of Bergman fans of matroids or polymatroids.
For simplicity, we only describe the Bergman fan of a loopless matroid. For more basics on matroids, see \cite{oxleyMatroidBook}.

Let $M=(E,r)$ be a matroid, where $E=[n]$ is the ground set and $r:2^{E}\rightarrow \mathbb{Z}_{\geq 0}$ is the rank function such that
\begin{itemize}
  \item $r(\emptyset)=0$;
  \item $r(I) \leq r(J)$, for any $I \subset J \subset E $;
  \item $r(I\cap J)+r(I\cup J) \leq r(I)+r(J)$, for any $I, J \subset E $.
\end{itemize}
The matroid $M$ is called loopless if we have further $r(i)=1$ for any $i \in E$.

A subset $F\subset E$ is called a flat if for any $F'\supsetneq F$, it holds that $r(F)<r(F')$. The collection of flats is denoted by $\mathcal{L}$. Denote $\mathcal{L}^*=\mathcal{L} \setminus \{\emptyset, E\}$ the set of proper flats.

Let $$\Delta_M ^*=\{\emptyset \subsetneq F_1 \subsetneq ...\subsetneq F_k \subsetneq E\}$$ be the set of flags of proper flats of $M$.
One can check by definition that $\Delta_M ^*$ is a simplicial complex of dimension $r(M)-2$, where $r(M)=r(E)$ is the rank of $M$.

The Bergman fan $\Sigma_M$ of $M$ a geometric realization of the complex $\Delta_M ^*$ defined as follows. For the ground set $E=[n]$, we fix a basis $\{e_i| i\in [n]\}$ of the vector space $\mathbb{R}^E$.
Put $V=\frac{\mathbb{R}^E}{\mathbb{R}{e_E}}$, and denote the image of $e_S$ by $u_S$ in the quotient space $V$. For a flag $\mathcal{F}=(F_1 \subsetneq ...\subsetneq F_k )\in \Delta_M ^*$, we associate a cone $$\sigma_{\mathcal{F}}=\cone(u_{F_1},...,u_{F_k}) \subset V.$$

Then $$\Sigma_M=\{\sigma_{\mathcal{F}}|\mathcal{F}\in \Delta_M ^*\}$$ is a simplicial fan of pure dimension $r(M)-1$ in $V$, called the Bergman fan of $M$. The Bergman fan is always balanced, hence, tropical.

\begin{prop}
The Bergman fan $\Sigma_M$ is Lorentzian.
\end{prop}

This follows from \cite{nowak2023mixed, ross2023lorentzianfan} or \cite{huhHRR}.
Furthermore, by \cite{huhHRR} the K\"ahler package holds on the Bergman fan $\Sigma_M$. This result was also extended to polymatroids, see \cite{pagaria2022hodgepolymatr} or \cite{huhwangPolymatroid}. In particular, the Bergman fan $\Sigma_M$ has hard Lefschetz property.

\section{Positivity criterion}\label{sec posicri}
The following result was proved in our previous work \cite{huxiaohardlef2022Arxiv}. Let $X$ be a smooth projective variety (or in the analytic setting, a compact K\"ahler manifold) of dimension $n$.

\begin{thrm} [Non-vanishing criterion] \label{nonvanishing}
  Let $m\leq n$ and let $\alpha_1,...,\alpha_m$ be nef classes on $X$. Then the following statements are equivalent:
\begin{enumerate}
  \item The intersection class $\alpha_1\cdot...\cdot\alpha_m \neq 0 $;
  \item For any subset $I \subset [m]$, one has $\nd(\alpha_I) \geq |I|$.
\end{enumerate}
\end{thrm}

Our proof for the above result applied the following proportionality statement, which has also appeared in \cite{luotieHodgeindex}(see also \cite{dinh_sibony_groupes_commutatifs, dangfarveBdivisor, zhang2020hrr}).

\begin{thrm}
[Proportionality]
\label{prop geom} Let $\mathfrak{L}=(L_1,...,L_{n-2})$ be a collection of nef classes on $X$. Assume that $\alpha, \gamma$ are nef class satisfying
$$(\mathfrak{L}\cdot\alpha\cdot \gamma)^2 = (\mathfrak{L}\cdot\alpha^2)(\mathfrak{L}\cdot \gamma^2).$$
Then $\mathfrak{L} \cdot \gamma$ and $\mathfrak{L} \cdot \alpha$ are proportional.
\end{thrm}

We extend these results to the setting of Lorentzian polynomials, which can be applied in various contexts.
The proofs are essentially the same as those in \cite{huxiaohardlef2022Arxiv}, for the reader's convenience, we include the details here.

\subsection{Proportionality}
We first give a proportionality result for Lorentzian quadratic forms (see \cite{luotieHodgeindex}).

\begin{lem}\label{luotie}
Let $(V,Q)$ be a finite-dimensional real vector space $V$ endowed with a quadratic form $Q$ having at most one positive eigenvalue. For $v\in V$, write $Q(v)=Q(v,v)$. Then the following results hold:
    \begin{enumerate}
        \item For any $v,w \in V$ with $Q(v)\geq 0$, we have $Q(v,w)^2\geq Q(v)Q(w)$.
        \item For any $v,w \in V$ with $Q(v)> 0$, the equality $Q(v,w)^2= Q(v)Q(w)$ holds if and only if $$w-\frac{Q(v,w)}{Q(v)}v\in \ker Q.$$
    \end{enumerate}
\end{lem}
\begin{proof}
(1) When $Q(v)=0$ or $v, w$ are proportional, this is clear.

We assume that $Q(v)>0$ and that $v, w$ are not proportional. Consider the function $q(t)=Q(w+tv)$. Since $Q$ has at most one positive eigenvalue, $Q$ cannot be nonnegative on the two-dimensional subspace spanned by $v, w$, hence $q(t)$ has at least one real root. This yields that
$$Q(v,w)^2\geq Q(v)Q(w).$$

(2) Write $x=w-\frac{Q(v,w)}{Q(v)}v$. Then the assumption $Q(v,w)^2= Q(v)Q(w)$ implies that $Q(v, x)=Q(x,x)=0$.

It suffices to show the following statement: whenever $Q(v, x)=Q(x,x)=0$, we have that $x\in \ker Q$.
Denote
$$ P=\ker Q(v, -) =\{y\in V| Q(v, y) =0\}.$$
By the assumption, $Q$ is positive definite on the line $\mathbb{R} v$. By (1), $Q$ is negative semi-definite on the hyperplane $P$. Therefore $Q$ has signature of the form $(+, -,...,-,0,...,0)$.

According to the signature, we fix a basis $e_k$ (with $e_1 = v$), $1\leq k \leq n$ such that $Q(e_k, e_l)=0$ for $k\neq l$, $Q(e_1, e_1)>0$, $Q(e_i, e_i)<0$ for $2 \leq i \leq s$ and $Q(e_i, -)=0$ for $s+1 \leq i \leq n$. Thus $Q(v, x)=0$ yields that
$$x= \sum_{ i\geq 2} c_i e_i.$$
If we further have $Q(x,x)= 0$, then
$$x= \sum_{i : \ Q(e_i, -)= 0} c_i e_i,$$
implying that $x\in \ker Q$.

This finishes the proof.
\end{proof}

\begin{prop}\label{proportional}
Let $(V,Q)$ be a finite-dimensional real vector space $V$ endowed with a quadratic form $Q$. Assume that $Q$ has exactly one positive eigenvalue. Assume that $x, y\in V$ satisfy $Q(x)\geq 0, Q(y)\geq 0$ and
$$Q(x, y)^2 = Q(x)Q(y),$$
then the linear forms $Q(x, -)$ and $Q(y, -)$ are proportional.
\end{prop}

\begin{proof}
We can assume that the linear forms $Q(x, -)$ and $Q(y, -)$ are non-zero, otherwise, it is obvious.

We first consider the case when at least one of $Q(x), Q(y)$ is nonzero.
Without loss of generality, assume that $Q(x)>0$. Note that
$$z= y- \frac{Q(x, y)}{Q(x)} x$$
satisfies $Q(z)=0, Q(x, z)=0$, then applying Lemma \ref{luotie} to $z$ implies:
$$ Q(y,-)=\frac{Q(x, y)}{Q(x)} Q(x,-).$$

Next, we consider the case when $Q(x)=Q(y)=0$. In particular, the equality $Q(x, y)^2 = Q(x)Q(y)$ implies that $Q(x, y)=0$. Consider the subspace spanned by $x,y$. If it is of dimension one, then there is nothing to prove. Therefore we may assume that it is of dimension two. Fix a vector $v$ such that $Q(v)>0$, then $\Span_{\mathbb{R}} \{x,y\}$ has a nontrivial intersection with the hyperplane $\ker Q(v,-)$, i.e., there exist two numbers $c, d$, not all zero, such that $Q(v, c x + dy)=0$. 
The aforementioned vanishing properties of $x, y$ imply that
\begin{equation*}
  Q(c x + dy)=0.
\end{equation*}
Then applying Lemma \ref{luotie} (2) (or its proof) shows that $Q(c x + dy, -)=0$, implying the result.

This finishes the proof.
\end{proof}

\subsection{Numerical dimension}

In the following, we fix an open convex cone $\mathfrak{C} \subset \mathbb{R}^s$ and a $\mathfrak{C}$-Lorentzian polynomial $f$ of degree $n$.
Let $F$ be the complete polarization of $f$, i.e., the symmetric multilinear function defined by
\begin{equation*}
  F(v_1,...,v_n)=\frac{1}{n!}D_{v_1}...D_{v_n}f.
\end{equation*}
We have that $ F(v_1,...,v_n) \geq 0$ when the $v_i \in \overline{\mathfrak{C}}$, where $\overline{\mathfrak{C}}$ is the closure of $\mathfrak{C}$.

Write $$F(v_1[n_1],...,v_d[n_d])=F(v_1,...,v_1,...,v_i,...,v_i,...,v_d,...,v_d)$$
when $v_i$ appears $n_i$ times and $\sum_i n_i=n$.

By the definition of Lorentzian polynomials, for any $v_1,...,v_n \in \overline{\mathfrak{C}}$,
\begin{equation*}
  F(v_1,v_2,v_3,...,v_n)^2 \geq F(v_1,v_1,v_3,...,v_n)F(v_2,v_2,v_3,...,v_n).
\end{equation*}
This is the Hodge index inequality (or Alexandrov-Fenchel inequality) for Lorentzian polynomials.

We introduce a notion of numerical dimension for $\mathfrak{C}$-Lorentzian polynomials.

\begin{lem}\label{DefnOfND}
Let $v\in \overline{\mathfrak{C}}$. Then the following number is independent of the choice of $w\in \mathfrak{C}$:
\begin{equation*}
  \nd_{w} (v)=\max \{k\in \mathbb{Z}_{\geq 0}: F(v[k],w[n-k])>0\}.
\end{equation*}
We denote this number by $\nd (v)$ and call it the numerical dimension of $v$.
\end{lem}

\begin{proof}
Let $w_1,w_2\in \mathfrak{C}$. It suffices to show $\nd_{w_2}(v)\geq \nd_{w_1}(v)=:d_0$, that is,  $$F(v[d_0],w_2[n-d_0])> 0.$$

Let $a_i=F(v[d_0],w_1[n-d_0-i],w_2[i]),i=0,...,n-d_0$. We only need to show that $a_{i}>0,$ for any $i$. The case $i=0$ follows from $\nd_{w_1}(v)=d_0$.
By induction on $i$, we may suppose $i>0$ and $a_{i-1}>0$.
  If $a_{i}=0$, then $a_{i-1}-ca_{i}>0$ for any $c\in \mathbb{R}$. On the other hand, note that
  \begin{equation*}
      a_{i-1}-ca_{i}=F(v[d_0],w_1[n-d_0-i],w_2[i-1],w_1-cw_2).
  \end{equation*}
  Since $\mathfrak{C}$ is open, we may take $c>0$ large enough such that $w_1-cw_2\in -\mathfrak{C}$. Then $a_{i-1}-ca_{i}\leq 0$, which gives a contradiction. So we must have $a_{i}>0$.

This finishes the proof.
\end{proof}

\begin{rmk}\label{submod}
In the setting when $F$ is given by the volume function of nef classes on a projective manifold $X$, $F$ is $\Amp^1(X)$-Lorentzian. The numerical dimension has an important submodularity property \cite{huxiaohardlef2022Arxiv}: for any $A, B, C\in \Nef^1 (X)$,
\begin{equation*}
  \nd(A+B+C)+\nd(C)\leq \nd(A+C)+\nd(B+C).
\end{equation*}
By the mechanism developed in \cite{huxiao2023intersectionLorentArxiv}, one checks that the same conclusion holds for a $\mathfrak{C}$-Lorentzian polynomial: for any $x, y, z\in \overline{\mathfrak{C}}$, we have that
\begin{equation*}
  \nd(x+y+z)+\nd(z)\leq \nd(x+z)+\nd(y+z).
\end{equation*}
\end{rmk}

\subsection{Non-vanishing criterion} \label{nonvanishing sec}

We prove the following result.

\begin{thrm}\label{PosiCri}
Let $m\in[n]$ and fix $\omega\in \mathfrak{C}$, then for any $\alpha_1,...,\alpha_m \in \overline{\mathfrak{C}}$,
\begin{equation*}
  F(\alpha_1,...,\alpha_m, \omega[n-m])>0 \Leftrightarrow \nd(\alpha_{I})\geq |I|,\forall I\subset [m].
\end{equation*}

\end{thrm}

\begin{proof}

It is clear that if $$F(\alpha_1,...,\alpha_m, \omega[n-m])>0,$$ then $\nd(\alpha_I) \geq |I|$ for any subset $I\subset [m]$.

For the converse direction, we prove the result by induction on $m$.

When $m=1$, this follows from the definition of $\nd(-)$.
Assume that the statement holds for $m-1$. In particular, we have
$$F(\alpha_1,...,\alpha_{m-1}, \omega[n-m+1])>0.$$

We first show that under the assumption on $I=[m]$, that is, $$F(\alpha_{[m]}[m], \omega [n-m])>0,$$
we must have 
\begin{equation}\label{eq m}
F(\alpha_1 , ..., \alpha_{m-1} , \alpha_{[m]}, \omega [n-m]) >0.
\end{equation}
Otherwise, assume that
\begin{equation}\label{eq otherwise}
  F(\alpha_1 , ..., \alpha_{m-1}, \alpha_{[m]} , \omega [n-m]) =0,
\end{equation}
then by applying Corollary \ref{proportional} to the quadratic form $$Q(-, -)=F(-,-,\alpha_2 , ..., \alpha_{m-1}, \omega [n-m])$$ with $x=\alpha_1, y=\alpha_{[m]}$ (by the Alexandrov-Fenchel inequality for $F$ and $Q(\omega)>0$, $Q$ has exactly one positive eigenvalue),
\begin{equation*}
F(\alpha_{[m]}- c \alpha_1, -, \alpha_2 , ..., \alpha_{m-1} , \omega [n-m]) = 0
\end{equation*}
for some real number $c$. It follows that $c>0$ by evaluating the above linear form at $\omega$ and using the induction.
 By (\ref{eq otherwise}), this implies that $$F(\alpha_{[m]}[2],\alpha_2 , ...,\alpha_{m-1} , \omega[n-m])=0. $$
Applying Corollary \ref{proportional} again to the quadratic form $F(-,-,\alpha_{[m]} ,\alpha_3, ..., \alpha_{m-1}, \omega [n-m])$ with $x=\alpha_2, y=\alpha_{[m]}$, we get
\begin{equation*}
 F(\alpha_{[m]} -c' \alpha_2,-,\alpha_{[m]} ,\alpha_3, ..., \alpha_{m-1}, \omega [n-m])=0
\end{equation*}
for some $c'>0$, yielding that
$$F(\alpha_{[m]} [3], \alpha_3 , ...,\alpha_{m-1} , \omega [n-m])=0.$$

By iteration, we obtain $$F (\alpha_{[m]}[m], \omega [n-m])=0,$$
which is a contradiction with the assumption on $I=[m]$.
Therefore, we must have
\begin{equation*}
  F(\alpha_1 , ..., \alpha_{m-1}, \alpha_{[m]},\omega [n-m]) >0.
\end{equation*}

We show that this implies $$F(\alpha_1 , ...,  \alpha_{m} , \omega [n-m]) >0,$$
which completes the induction.
Otherwise,
\begin{equation}\label{eq otherwise2}
  F(\alpha_1 , ...,  \alpha_{m} , \omega [n-m])=0.
\end{equation}
Then by Corollary \ref{proportional}, we get that the forms $$F(\alpha_1,-, \alpha_3, ...,\alpha_{m}, \omega [n-m]),\ F(\alpha_2,-, \alpha_3, ...,\alpha_{m}, \omega [n-m])$$
are proportional. By inductive assumption on $m-1$, we have:
$$F(\alpha_1,\omega, \alpha_3, ...,\alpha_{m}, \omega [n-m])>0$$ and
$$F(\alpha_2,\omega, \alpha_3, ...,\alpha_{m}, \omega [n-m])>0,$$
therefore, there is some $c''>0$ such that
\begin{equation*}
F(\alpha_1 -c'' \alpha_2,-, \alpha_3, ...,\alpha_{m}, \omega [n-m])=0,
\end{equation*}
which by (\ref{eq otherwise2}) yields that $$F(\alpha_1 [2] ,\alpha_3, ..., \alpha_{m} , \omega [n-m]) =0.$$

The same argument shows that
\begin{equation*}
 F(\alpha_1, ..., \alpha_{m-1}, \alpha_{k}, \omega [n-m]) =0
\end{equation*}
for every $1\leq k \leq m-1$.
Together with (\ref{eq otherwise2}), we get that
$$F(\alpha_1 , ..., \alpha_{m-1} , \alpha_{[m]} , \omega [n-m]) =0,$$
which contradicts (\ref{eq m}).

This finishes the proof.
\end{proof}

\begin{exmple}
By specifying the $\mathfrak{C}$-Lorentzian polynomial $f$ in different settings, we obtain a non-vanishing criterion for complete intersections in various contexts. For example, by taking $f$ to be the volume function with $\mathfrak{C}$ given by the ample or K\"ahler cone, it holds on smooth projective varieties, compact K\"ahler manifolds and Lorentzian fans.
\end{exmple}

\section{Basics on collections of nef classes}\label{sec collec}

In this section, we present some basic properties of collections of nef classes, some of which will be applied in Section \ref{sec arbi}. 

While the results hold in various contexts beyond geometry, the reader may keep in mind the algebro-geometric situation. In this section, we fix a smooth projective variety $X$ of dimension $n$ and assume that the collection $\mathfrak{L}=(L_1,...,L_{m})$ is a subset of $\Nef^1 (X)$ with $m\leq n$. We establish several statements analogous to those in \cite{handel2022AFextremals}.

Recall the definition from the introduction:

\begin{defn}
A collection $\mathfrak{L}=(L_1,...,L_{m})$ of nef classes is called
\begin{itemize}
  \item subcritical, if $\nd(L_I) \geq |I|$ for any $I \subset [m]$;
  \item critical, if $\nd(L_I) \geq |I|+1$ for any $I \subset [m]$;
  \item supercritical, if $\nd(L_I) \geq |I|+2$ for any $I \subset [m]$.
\end{itemize}
\end{defn}

It is clear that: supercritical $\Rightarrow$ critical $\Rightarrow$ subcritical.

By using Theorem \ref{nonvanishing} or \ref{PosiCri}, the collection $\mathfrak{L}$ is subcritical if and only if the complete intersection class $$\mathfrak{L}=L_1 \cdot... \cdot L_{m}\neq 0.$$

We also have:

\begin{prop}\label{kernel nontriv}
Let $\mathfrak{L}=(L_1,...,L_{n-2})$ be a subcritical collection of nef classes on $X$. Assume that the collection $\mathfrak{L}$ is not critical. Then $\ker \mathfrak{L}$ is non-trivial. In particular, the class $\mathfrak{L}$ can be an HL class only if the collection $\mathfrak{L}$ is at least critical.
\end{prop}

\begin{proof}
By the assumption, there is a non-empty subset $I\subset [n-2]$ such that $\nd(L_I) =|I|$. By Theorem \ref{nonvanishing} or \ref{PosiCri}, this yields $\mathfrak{L} \cdot L_i =0$ for any $i\in I$. Therefore, $\{L_i\}_{i\in I} \subset \ker \mathfrak{L}$.
\end{proof}

\subsection{The critical case}
Throughout this subsection, we assume that $\mathfrak{L}=(L_1,...,L_m)$ is a critical but not supercritical collection of nef classes. This means that $\nd(L_I)\geq |I|+1$ for all $I\subset [m]$ and $\nd(L_{I_0})=|I_0|+1$ for some $I_0 \subset [m]$.

\begin{defn}
A subset $I\subset [m]$ is called $\mathfrak{L}$-critical if $\nd(L_I)=|I|+1$, and $I$ is called $\mathfrak{L}$-maximal if it is $\mathfrak{L}$-critical and there is no $\mathfrak{L}$-critical subset $J$ such that $I \subsetneq J$.
\end{defn}

\begin{lem}
If $I, J$ are $\mathfrak{L}$-critical and $K=I\cap J \neq\emptyset$, then $I\cup J$ is $\mathfrak{L}$-critical.
\end{lem}

\begin{proof}
Write $I, J$ as disjoint unions: $I =I' \sqcup K$, $J=J' \sqcup K$, then $I\cup J =I'\sqcup J'\sqcup K$. Since $\mathfrak{L}$ is a critical collection, we have $$\nd(L_{I\cup J})\geq |I\cup J|+1 =|I'|+ |J'|+|K|+1.$$
On the other hand, by using the inequality $\nd(L_K)\geq |K| +1$ and the submodularity of numerical dimension \cite{huxiaohardlef2022Arxiv} (see Remark \ref{submod}),
\begin{align*}
  \nd(L_{I\cup J}) &\leq \nd(L_I) +\nd(L_J) -\nd(L_K) \\
&= |I\cup J|+1 + |K|+1 - \nd(L_K) \\
&\leq |I\cup J|+1.
\end{align*}
Thus, $I\cup J$ is $\mathfrak{L}$-critical.
\end{proof}

\begin{cor}
If $I $ and $ J$ are distinct $\mathfrak{L}$-maximal, then $I\cap J=\emptyset$ and $\nd(L_{I\cup J})=|I|+|J|+2$.
\end{cor}

\begin{proof}
If $I\cap J\neq \emptyset$, then by the above lemma, $I\cup J$ must be $\mathfrak{L}$-critical, contradicting with the $\mathfrak{L}$-maximality of $I, J$.
Therefore, $I\cap J=\emptyset$. Again, by the $\mathfrak{L}$-maximality of $I$ and $J$, $I\cup J$ cannot be $\mathfrak{L}$-critical, yielding that $$\nd(L_{I\cup J})\geq |I\cup J|+2=|I|+|J| +2.$$
Finally, note that $$\nd(L_{I\cup J}) \leq \nd(L_I) +\nd(L_J) = |I|+|J|+2,$$ since $I, J$ are $\mathfrak{L}$-critical.

This finishes the proof.
\end{proof}

\begin{cor}
For a critical but non-supercritical collection $\mathfrak{L}=(L_1,...,L_m)$, there exist nonempty $\mathfrak{L}$-maximal subsets, and any two different $\mathfrak{L}$-maximal subsets are disjoint.
\end{cor}


\subsection{The subcritical case}
Throughout this subsection, we assume that the collection of nef classes $\mathfrak{L}=(L_1,...,L_m)$ is subcritical but not critical, that is, $\nd(L_I)\geq |I|$ for all $I\subset [m]$ and $\nd(L_{I_0})=|I_0|$ for some nonempty subset $I_0 \subset [m]$.

\begin{defn}
A subset $I\subset [m]$ is called $\mathfrak{L}$-subcritical if $\nd(L_I)=|I|$.
\end{defn}

\begin{lem}\label{subcr union}
If $I, J$ are $\mathfrak{L}$-subcritical, then $I\cup J$ is also $\mathfrak{L}$-subcritical.
\end{lem}

\begin{proof}
By the criticality of $\mathfrak{L}$,
\begin{equation*}
  \nd(L_{I\cup J})\geq |I\cup J|
\end{equation*}

On the other hand, by using the inequality $\nd(L_{I\cap J})\geq |I\cap J|$ and the submodularity of numerical dimension (see Remark \ref{submod}),
\begin{equation*}
  \nd(L_{I\cup J}) \leq \nd(L_I) +\nd(L_J) -\nd(L_{I\cap J}) \leq |I|+|J| - |I\cap J|.
\end{equation*}
Thus, $\nd(L_{I\cup J})= |I\cup J|$, i.e., $I\cup J$ is $\mathfrak{L}$-subcritical.
\end{proof}

\begin{cor}\label{max subcritical}
There exists a unique maximal subcritical subset $M$, which contains every $\mathfrak{L}$-subcritical subset.
\end{cor}

\begin{proof}
By Lemma \ref{subcr union}, we can take $M$ to be the union of all the subcritical subsets of $\mathfrak{L}$.
\end{proof}

\begin{prop}\label{subcr to cri}
Let $M$ be the maximal subcritical subset. Then the collection $\mathfrak{L}_{\setminus M}$ is critical on the class $\mathfrak{L}_M$, in the sense that, for any $I\subset [n-2]\setminus M$,
\begin{equation*}
  L_I ^{|I|+1} \cdot \mathfrak{L}_M \neq 0.
\end{equation*}
\end{prop}

\begin{proof}
We argue by contradiction. Suppose there exists $I\subset [n-2]\setminus M$ such that
\begin{equation}\label{subcr vani}
  L_I ^{|I|+1} \cdot \mathfrak{L}_M = 0.
\end{equation}

Denote $Q$ the collection of nef classes given by $L_I$ repeating $|I|+1$ times and $L_m$, $m\in M$. By Theorem \ref{nonvanishing} or \ref{PosiCri}, (\ref{subcr vani}) holds if and only if there exists some subcollection $K\subset Q$ such that
\begin{equation*}
  \nd(\sum_{k\in K} L_k)< |K|.
\end{equation*}
An easy analysis shows that this is possible only if $K\subset Q$ is given by $L_I$ repeating $|I|+1$ times and $L_{m_1},...,L_{m_k}$ with $m_i \in M$, therefore,
\begin{equation*}
  \nd(\sum_{k\in K} L_k)=\nd(L_I + L_{m_1}+...+L_{m_k})< |I|+1+k.
\end{equation*}
On the other hand, by the criticality of $\mathfrak{L}$,
\begin{equation*}
\nd(L_I + L_{m_1}+...+L_{m_k})\geq |I|+k.
\end{equation*}
This yields that $I\cup \{m_1,...,m_k\}$ must be a subcritical set. This contradicts the maximal property of $M$ -- that is, $M$ contains all the  subcritical sets.

This finishes the proof.
\end{proof}

\section{The local Hodge index inequality}\label{sec localaf}

\subsection{General form for Lorentzian polynomials}

Given a finite-dimensional real vector space $V$ and a subset $S \subset V$, we use the following notation:
\begin{itemize}
  \item $\overline{{\cone}}(S)$ is the closed convex cone generated by $S$, that is, the closure of the set of convex combinations of elements in $S$.
\end{itemize}
The dual space of $V$ is denoted by $V^*$. Given a convex cone $\mathfrak{C}\subset V$, its dual cone, denoted by $\mathfrak{C}^*$, is the set of linear functionals that take nonnegative values on $\mathfrak{C}$.

In this section, we first prove a local Hodge index inequality in a general form, inspired by \cite[Theorem 4.3]{handel2022AFextremals}.
To simply the notation and present the results in a more geometric form, we introduce the following terminology. Let $\mathfrak{C}$ be an open convex cone in $V=\mathbb{R}^s$, and let $f$ be a $\mathfrak{C}$-Lorentzian polynomial of degree $n$ on $V$. Let $F$ be the complete polarization of $f$. Given $\alpha_1,...,\alpha_n \in V$, we write
\begin{equation*}
  F(\alpha_1,...,\alpha_n)=\alpha_1 \cdot ...\cdot \alpha_n.
\end{equation*}

\subsubsection{Positive cones for Lorentzian polynomials} \label{pos cone lor}

Analogous to the geometric setting, we set
$$\Amp^1 (f)=\mathfrak{C},\ \Nef^1 (f)=\overline{\mathfrak{C}}.$$
Elements in $\Amp^1 (f)$ (resp. $\Nef^1 (f)$) are called ample/K\"ahler (resp. nef) classes with respect to $f$, or simply ample/K\"ahler (resp. nef) classes.  We let $\Eff^1(f) \subset V$ be a closed convex cone such that $\Nef^1 (f) \subset \Eff^1(f)$. We further require that
$$T\cdot\alpha_1\cdot ...\cdot \alpha_{n-1}\geq 0,$$
for any $T\in \Eff^1(f)$ and $\alpha_1, ..., \alpha_{n-1}\in \Nef^1 (f)$. Elements in $\Eff^1(f)$ are called pseudo-effective classes, and the interior points of $\Eff^1(f)$ are called big classes.

\begin{rmk}

Using Lemma \ref{DefnOfND}, the definitions of supercritical, critical, and subcritical collections of nef classes in the geometric setting extend naturally to Lorentzian polynomials.

\end{rmk}

Now we can state our result:

\begin{thrm}\label{local af}
Let $\mathfrak{L}=(L_1,...,L_{n-2})$ be a critical collection of nef classes in $\Nef^1 (f)$. Let $\mathcal{C} \subset \Eff^1 (f) $ be a collection of pseudo-effective classes satisfying that
\begin{itemize}
  \item $\overline{\cone} (\mathcal{C} ) = \Eff^1 (f) $;
  \item $ \mathcal{C} \cap \ker \mathfrak{L}$ is a finite set;
  \item for any $D\in \mathcal{C}$ and $\alpha_1,...,\alpha_{n-3} \in \Nef^1 (f)$, the product $D\cdot \alpha_1 \cdot... \cdot\alpha_{n-3}$ has the weak Hodge-Riemann property, that is,
      the quadratic form
      $$q(x,y) = x\cdot y \cdot D\cdot \alpha_1 \cdot... \cdot\alpha_{n-3}$$
      has at most one positive eigenvalue.
\end{itemize}
Fix $r\in [n-2]$. Then for any $\alpha \in \ker \mathfrak{L}$, there exists $\beta \in V$ such that
\begin{enumerate}
  \item $\beta - \alpha \in  \Span_\mathbb{R} (\mathcal{C}\cap \ker \mathfrak{L})$, and
  \item $-\beta^2 \cdot \mathfrak{L}_{\setminus r} \in \Eff^1 (f)^*$.
\end{enumerate}

\end{thrm}

By $\alpha \in \ker \mathfrak{L}$ we mean that $l(-)=\mathfrak{L} \cdot \alpha\cdot - =0$ as an element of $V^*$.

\begin{proof}
Let $A$ be a fixed ample or K\"ahler class.

\begin{claim}
There exists $\beta$ such that $\beta - \alpha \in  \Span_\mathbb{R} (\mathcal{C}\cap \ker \mathfrak{L})$ and
\begin{equation}\label{af matrix}
  D\cdot \beta \cdot A\cdot \mathfrak{L}_{\setminus r} =0
\end{equation}
holds for any $D\in \mathcal{C} \cap \ker \mathfrak{L}$.
\end{claim}

Assuming the claim, we show that the class $\beta$ constructed in the claim satisfies that $$-\beta^2 \cdot \mathfrak{L}_{\setminus r} \in \Eff^1 (f)^*,$$
which, by the assumption $\overline{\cone} (\mathcal{C} ) = \Eff^1 (f)$, is equivalent to
\begin{equation}\label{movable}
  D \cdot \beta^2 \cdot \mathfrak{L}_{\setminus r} \leq 0,\ \forall D \in \mathcal{C}.
\end{equation}

To prove (\ref{movable}), we consider the following two alternatives on $D \in \mathcal{C}$.

\begin{description}
  \item [Case 1] $\mathfrak{L} \cdot D =0$.
\end{description}

By (\ref{af matrix}) and the assumption that $D\cdot \mathfrak{L}_{\setminus r}$ has weak Hodge-Riemann property, we get that
\begin{equation*}
  0=\left(D\cdot \beta \cdot A\cdot \mathfrak{L}_{\setminus r}\right)^2 \geq \left(D\cdot \beta^2 \cdot \mathfrak{L}_{\setminus r}\right)\left(D\cdot A^2 \cdot \mathfrak{L}_{\setminus r}\right).
\end{equation*}
If $D \cdot A^2 \cdot \mathfrak{L}_{\setminus r} >0$, then it is clear that
$$D\cdot \beta^2 \cdot \mathfrak{L}_{\setminus r} \leq 0,$$
as desired.
If $D\cdot A^2 \cdot \mathfrak{L}_{\setminus r} =0$, then we must have
$$ D \cdot A \cdot \mathfrak{L}_{\setminus r} =0.$$
This is because $D \cdot A \cdot \mathfrak{L}_{\setminus r} \in \Nef^1(f)^* \subset V^*$ and $A \in \Amp^1(f)$.
It follows that $$ D \cdot \beta ^2 \cdot \mathfrak{L}_{\setminus r}=D \cdot (\beta + tA)^2 \cdot \mathfrak{L}_{\setminus r}$$ for any $t\in \mathbb{R}$.
Since $\Amp^1(f)$ is open, we can choose $t>0$ large enough such that $\widetilde{A}=\beta +tA$ is ample.
Finally, note that for any $A,\widetilde{A}\in \Amp^1(f)$,
\begin{equation*}
  D\cdot A^2 \cdot \mathfrak{L}_{\setminus r} =0 \iff D\cdot \widetilde{A}^2 \cdot \mathfrak{L}_{\setminus r} =0.
\end{equation*}

In turn, we find
$$D\cdot \beta^2 \cdot \mathfrak{L}_{\setminus r} = D\cdot (\beta +tA)^2 \cdot \mathfrak{L}_{\setminus r}=0.$$

This finishes the proof of Case 1.

\begin{description}
  \item[Case 2] $\mathfrak{L} \cdot D \neq0$.
\end{description}

Since $\beta - \alpha \in  \Span_\mathbb{R} (\mathcal{C}\cap \ker \mathfrak{L})$ and $\alpha \in \ker \mathfrak{L}$, we must have $\mathfrak{L} \cdot \beta =0$. Then using again the weak Hodge-Riemann property of $D\cdot \mathfrak{L}_{\setminus r}$, we get that
\begin{equation*}
  0=\left(D\cdot \beta \cdot L_r \cdot \mathfrak{L}_{\setminus r}\right)^2 \geq  \left(D\cdot \beta^2 \cdot \mathfrak{L}_{\setminus r}\right)\left(D\cdot L_r ^2\cdot \mathfrak{L}_{\setminus r}\right).
\end{equation*}

If $D\cdot L_r ^2\cdot \mathfrak{L}_{\setminus r}>0$, then it is clear that
$$D\cdot \beta^2 \cdot \mathfrak{L}_{\setminus r}\leq 0.$$

It remains to consider the case when $D \cdot L_r ^2\cdot \mathfrak{L}_{\setminus r}=0$. Since $\mathfrak{L} \cdot D \neq0$,
\begin{equation}\label{positive}
  D\cdot A \cdot \mathfrak{L} = D\cdot A \cdot L_r \cdot \mathfrak{L}_{\setminus r}>0,
\end{equation}
which yields that there is some constant $b\in \mathbb{R}$ such that
\begin{equation}\label{af b}
  D\cdot (\beta -bL_r) \cdot A \cdot \mathfrak{L}_{\setminus r}=0.
\end{equation}
By (\ref{af b}) and the weak Hodge-Riemann property of $D\cdot \mathfrak{L}_{\setminus r}$,
\begin{equation}\label{af b1}
  0 =\left(D\cdot (\beta -bL_r) \cdot A \cdot \mathfrak{L}_{\setminus r}\right)^2\geq \left(D\cdot (\beta -bL_r)^2 \cdot \mathfrak{L}_{\setminus r}\right)\left(D\cdot A ^2\cdot \mathfrak{L}_{\setminus r}\right).
\end{equation}
Note that there is some constant $c>0$ such that $A - cL_r$ is nef, therefore by (\ref{positive}) we get
\begin{equation}\label{postive a}
 D\cdot A ^2\cdot \mathfrak{L}_{\setminus r} \geq c D\cdot A \cdot L_r\cdot \mathfrak{L}_{\setminus r} >0.
\end{equation}
Applying (\ref{af b1}), (\ref{postive a}), $\mathfrak{L}\cdot \beta =0$ and the assumption $D\cdot L_r ^2\cdot \mathfrak{L}_{\setminus r}=0$, we get that
\begin{equation*}
  D\cdot \beta ^2 \cdot \mathfrak{L}_{\setminus r}=D\cdot (\beta -bL_r)^2 \cdot \mathfrak{L}_{\setminus r} \leq 0.
\end{equation*}

This finishes the proof of Case 2.

Next we give the proof of the Claim.
\begin{proof}[Proof of the Claim]
By the assumption on $\mathcal{C}$, there is a finite collection $$\{D_i\}_{i\in [N]} \subset \mathcal{C}$$
that spans the vector space $V$ and contains all elements of $\mathcal{C} \cap \ker \mathfrak{L}$.

We fix such a collection. Then for any $\alpha \in \ker \mathfrak{L}$, there are some $a_i \in \mathbb{R}$ such that
\begin{equation}\label{vect alpha}
  \alpha =\sum_{i\in [N]} a_i D_i.
\end{equation}
We divide $[N]$ into two disjoint parts $[N]= V \cup V^c$, where
\begin{align*}
  &V=\{i\in [N]: \mathfrak{L} \cdot D_i \neq 0\},\\
 &V^c=\{i\in [N]: \mathfrak{L} \cdot D_i = 0\}=\mathcal{C} \cap \ker \mathfrak{L}.
\end{align*}

Let $\mathbb{R}^{[N]}$ be the real vector space of dimension $N$ by viewing $\{D_i: i \in [N]\}$ as a basis. Similarly, we have the vector spaces $\mathbb{R}^{V}$, $\mathbb{R}^{V^c}$, which are subspaces of $\mathbb{R}^{[N]}$. We endow $\mathbb{R}^{[N]}$ with the standard Euclidean inner product $\langle -, -\rangle$. Let $$p_V : \mathbb{R}^{[N]}\rightarrow \mathbb{R}^{V},\ p_{V^c} : \mathbb{R}^{[N]}\rightarrow \mathbb{R}^{V^c}$$ be the projections.

Let $\mathcal{A}$ be the matrix given by
\begin{equation}\label{matrix a}
  \mathcal{A}=[D_i \cdot D_j \cdot A \cdot \mathfrak{L}_{\setminus r}]_{i, j \in [N]}.
\end{equation}

Given $\alpha \in \ker \mathfrak{L}$, we fix $\mathbf{a}=(a_1,...,a_N)^t\in \mathbb{R}^{[N]}$ -- a column vector given by (\ref{vect alpha}). Then it is easy to see that the existence of $$\beta = \sum_{i\in[N]} b_i D_i$$ with the desired requirement in the Claim is equivalent to the solvability of the system of equations:
\begin{equation}\label{eq z}
  p_V (\mathbf{z})=p_V (\mathbf{a}),\ \text{and}\ p_{V^c} \mathcal{A} (\mathbf{z})=0.
\end{equation}

Next we show that (\ref{eq z}) is always solvable. Write $\mathbf{z} = p_V (\mathbf{z}) + p_{V^c} (\mathbf{z})$, it is sufficient to show that the equation
\begin{equation}\label{eq z2}
 p_{V^c} \mathcal{A} p_{V^c} (\mathbf{z}) =- p_{V^c} \mathcal{A} p_{V} (\mathbf{a})
\end{equation}
is solvable, which is equivalent to
\begin{equation}\label{eq z3}
  \image (p_{V^c} \mathcal{A} p_{V}) \subset \image (p_{V^c} \mathcal{A} p_{V^c}).
\end{equation}
Since the projection operators $p_V, p_{V^c}$ and the operator $\mathcal{A}$ are all self-adjoint, by taking orthogonal complements, (\ref{eq z3}) is equivalent to
\begin{equation}\label{eq z4}
  \ker (p_{V^c} \mathcal{A} p_{V^c}) \subset (\image (p_{V^c} \mathcal{A} p_{V}))^\perp.
\end{equation}

To prove the above inclusion, consider the subspace $$W=\{\mathbf{w}\in \mathbb{R}^{[N]}: p_{V^c} (\mathbf{w}) \in \ker \mathcal{A}\}.$$
Then it is clear that $W  \subset (\image (p_{V^c} \mathcal{A} p_{V}))^\perp$. We shall show that $$\ker (p_{V^c} \mathcal{A} p_{V^c}) \subset W.$$

Take $\mathbf{z}\in \ker (p_{V^c} \mathcal{A} p_{V^c})$, then
\begin{equation}\label{eq z5}
  \langle \mathcal{A} p_{V^c}(\mathbf{z}), p_{V^c} (\mathbf{z})\rangle=0
\end{equation}
Let $\gamma$ be the class given by $p_{V^c}(\mathbf{z}) \in \mathbb{R}^{V^c}$, then by the definition (\ref{matrix a}) of $\mathcal{A}$, (\ref{eq z5}) is equivalent to 
\begin{equation}\label{eq gamma}
  \gamma^2 \cdot A\cdot \mathfrak{L}_{\setminus r}=0.
\end{equation}
On the other hand, since $\gamma \in \mathbb{R}^{V^c}$, we have $\mathfrak{L} \cdot \gamma=0$, which in turn yields that
\begin{equation}\label{eq gamma1}
  \gamma \cdot L_r\cdot A\cdot \mathfrak{L}_{\setminus r}=0.
\end{equation}
Recall that $\mathfrak{L}$ is critical, thus
\begin{equation}\label{eq lr}
  L_r ^2 \cdot A\cdot \mathfrak{L}_{\setminus r} >0.
\end{equation}

By applying Proposition \ref{proportional} to (\ref{eq gamma}), (\ref{eq gamma1}) and (\ref{eq lr}), we get that
\begin{equation*}
  \gamma\cdot A \cdot \mathfrak{L}_{\setminus r} =0
\end{equation*}
which is equivalent to $p_{V^c} (\mathbf{z}) \in \ker \mathcal{A}$. This implies (\ref{eq z4}) and finishes the proof of the Claim.
\end{proof}

This completes the proof of the theorem.

\end{proof}

\subsection{Existence of the generating set}\label{sec generating set}

We call the subset $\mathcal{C} \subset \Eff^1 (f)$ a generating set when it satisfies the three requirements in Theorem \ref{local af}:
\begin{enumerate}
  \item $\overline{\cone} (\mathcal{C} ) = \Eff^1 (f) $;
  \item $ \mathcal{C} \cap \ker \mathfrak{L}$ is a finite set;
  \item for any $D\in \mathcal{C}$ and $\alpha_1,...,\alpha_{n-3} \in \Nef^1 (f)$, the product $D\cdot \alpha_1 \cdot... \cdot\alpha_{n-3}$ has weak Hodge-Riemann property.
\end{enumerate}
In this section, we establish the existence of such a generating set in various contexts. The existence and explicit description of such $\mathcal{C}$ are crucial in the numerical characterization of HL classes.

\subsubsection{Divisor classes on projective varieties}
Recall that we work on smooth projective varieties $X$ over an algebraically closed field of characteristic 0. 
The function $f=\vol$ is given by the volume function of divisor classes, which is defined on the vector space $V=N^1 (X)$. We let
\begin{itemize}
  \item $\mathfrak{C}=\Amp^1 (X)$,
  \item $\Eff^1 (f)=\Eff^1 (X)$,
\end{itemize}
as introduced in Section \ref{sec pos}. Then it is clear that $f$ is $\Amp^1 (X)$-Lorentzian, and the cones satisfy the condition required in Section \ref{pos cone lor}.

\begin{claim}
Let $\prim^1 (X)$ be the set of prime divisor classes, then $\prim^1 (X)$ is a generating set.
\end{claim}

Consider the three requirement in the definition of a generating set. The first is simply the definition of the pseudo-effective cone. The third follows from the Hodge index theorem. It remains to show that $\prim^1 (X) \cap \ker \mathfrak{L}$ is finite. 
This is established later in Proposition \ref{finite eff}.

\begin{lem}\label{lem mov restr}
Let $A$ be a very ample line bundle on $X$, and let $\alpha \in \Mov^1 (X)$. Then for a very general smooth hypersurface $H\in |A|$, the restriction $\alpha_{|H} \in \Mov^1 (H)$.
\end{lem}

\begin{proof}
This follows from the definition of movable divisor classes.
\end{proof}

\begin{prop}\label{kernel mov}
For any critical collection of nef classes $\mathfrak{L}=(L_1,...,L_{n-2})$, we have
$$\Mov^1 (X)\cap \ker \mathfrak{L}=\{0\}.$$
\end{prop}

We conjecture that the same statement holds in positive characteristic. However, we are currently unable to prove this due to the absence of a Lefschetz hyperplane theorem for $N^1(X)$ in this setting.

\begin{proof}
Fix a very ample line bundle $A$ on $X$. Let $\alpha \in \Mov^1 (X)$ satisfy 
\begin{equation}\label{movable a}
  \mathfrak{L}\cdot \alpha=0.
\end{equation}

We prove the result by induction on $n$.

When $n=3$, by the criticality assumption, $\mathfrak{L}$ is given by a nef class $L$ with $\nd(L) \geq 2$. By restricting (\ref{movable a}) to a very general section $H\in |A|$, we obtain $$L_{|H}\cdot \alpha_{|H} =0.$$
By Lemma \ref{lem mov restr}, $\alpha_{|H}$ is movable on the surface $H$. As on projective surfaces movable classes are nef, $\alpha_{|H}$ is also nef on $H$. Since $\nd(L) \geq 2$, the restriction $L_{|H}$ is big on $H$. By the duality of cones, we get that $\alpha_{|H} =0$. The Lefschetz hyperplane theorem \cite{fernexlauLefHyper, grothendieckLef} then implies that $\alpha =0$ on $X$.

When $n\geq 4$, assume that the result holds in $\dim =n-1$.

Fix $r \in [n-2]$. By restricting to a very general section $H\in |A|$, we know that $\alpha_{|H} \in \Mov^1 (H)$. On one hand, we have that
\begin{equation*}
  \alpha_{|H} \cdot {L_r}_{|H} \cdot {\mathfrak{L}_{\setminus r} }_{|H} =0,\  {L_r}_{|H} ^2 \cdot {\mathfrak{L}_{\setminus r} }_{|H} >0,
\end{equation*}
where the latter follows since $\mathfrak{L}$ is critical. Note that ${\mathfrak{L}_{\setminus r} }_{|H}$ is a product of nef classes on $H$, thus it has the weak Hodge-Riemann property, in particular,
\begin{equation*}
  0=\left( \alpha_{|H} \cdot {L_r}_{|H} \cdot {\mathfrak{L}_{\setminus r} }_{|H} \right)^2 \geq \left( \alpha_{|H}^2 \cdot {\mathfrak{L}_{\setminus r} }_{|H} \right)\left({L_r}_{|H} ^2 \cdot {\mathfrak{L}_{\setminus r} }_{|H} \right).
\end{equation*}
This yields that
\begin{equation*}
  \alpha_{|H} ^2 \cdot  {\mathfrak{L}_{\setminus r} }_{|H} \leq 0.
\end{equation*}
On the other hand, $\alpha_{|H} \in \Mov^1 (H)$ implies that $\alpha_{|H}  \cdot  {\mathfrak{L}_{\setminus r} }_{|H}$ has non-negative intersection with any element of $\Eff^1 (H)$, thus
$$\alpha_{|H}  \cdot  {\mathfrak{L}_{\setminus r} }_{|H} \in \Mov_1 (H). $$ Therefore,
\begin{equation*}
  \alpha_{|H} ^2 \cdot  {\mathfrak{L}_{\setminus r} }_{|H} \geq 0.
\end{equation*}
In summary, we have:
\begin{itemize}
  \item $\alpha_{|H} \cdot {L_r}_{|H} \cdot {\mathfrak{L}_{\setminus r} }_{|H} =0$;
  \item $\alpha_{|H} ^2 \cdot  {\mathfrak{L}_{\setminus r} }_{|H} = 0$;
  \item ${L_r}_{|H} ^2 \cdot {\mathfrak{L}_{\setminus r} }_{|H} >0$.
\end{itemize}
By Theorem \ref{prop geom}, there exists a constant $c\in \mathbb{R}$ such that 
$$\alpha_{|H} \cdot {\mathfrak{L}_{\setminus r} }_{|H}=c {L_r}_{|H} \cdot {\mathfrak{L}_{\setminus r} }_{|H}.$$
Multiplying both sides by $L_{r|H}$ implies $c=0$. Hence,
$$\alpha_{|H} \cdot  {\mathfrak{L}_{\setminus r} }_{|H} =0.$$

Since the collection ${\mathfrak{L}_{\setminus r} }_{|H}$ is critical on $H$, by induction, $\alpha_{|H}=0$. The Lefschetz hyperplane theorem then implies $\alpha =0$. This finishes the proof.

\end{proof}

As a consequence of Proposition \ref{kernel mov}, we obtain:

\begin{cor}\label{mov pos}
Let $\mathfrak{L}=(L_1,...,L_{n-2})$ be a critical collection of nef classes on $X$. Then for any big movable class $B$ on $X$, there is a unique big and movable class $M$ such that
\begin{equation*}
  B\cdot \mathfrak{L} = \langle M^{n-1}\rangle,
\end{equation*}
where $\langle - \rangle$ is the positive product of pseudo-effective classes.
\end{cor}

For the theory of positive products of pseudo-effective classes, we refer the reader to \cite{BDPP13, BFJ09}.

\begin{proof}
It is clear that $B\cdot \mathfrak{L} \in \Mov_1 (X)$.
By \cite[Section 3]{lehmannXiaoPosiCurve}, to show the existence of $M$, it is sufficient to show that any $\alpha\in \Mov^1 (X)$ satisfying $$B\cdot \mathfrak{L}\cdot \alpha =0$$ is zero.

Since $\mathfrak{L}\cdot \alpha \in \Mov_1 (X)$ and $B$ is big, by the cone duality $\Eff^1 (X)^* = \Mov_1 (X)$, $B\cdot \mathfrak{L}\cdot \alpha =0$ implies that $\mathfrak{L}\cdot \alpha =0$. By Proposition \ref{kernel mov}, $\alpha=0$.

The uniqueness of such $M$ also follows from \cite{lehmannXiaoPosiCurve}.
\end{proof}

\begin{prop}\label{finite eff}
Let $\mathfrak{L}=(L_1,...,L_{n-2})$ be a critical collection of nef classes on $X$. Then the set of prime divisors $D$ such that $\mathfrak{L} \cdot [D] =0$ is finite.
\end{prop}

\begin{proof}
By Corollary \ref{mov pos}, any prime divisor $D$ such that $\mathfrak{L} \cdot [D] =0$ satisfies that
\begin{equation*}
  \langle M^{n-1}\rangle \cdot [D] =0,
\end{equation*}
where $M$ is the class given in Corollary \ref{mov pos}.

By \cite{BFJ09}, this is equivalent to $D$ being a divisorial component of $\mathbb{B}_+ (M)$, where $\mathbb{B}_+ (M)$ is the augmented base locus of $M$. Therefore, the collection of such prime divisors is finite.
\end{proof}

\begin{rmk}
By \cite{Nak04} or \cite{Bou04}, the prime divisors in Proposition \ref{finite eff} are linearly independent and span extremal rays of the pseudo-effective cone. In particular, their number does not exceed the Picard number $\rho (X)$ of $X$.
\end{rmk}

\subsubsection{Lorentzian fans}

Let $\Sigma$ be a Lorentzian $n$-fan. In this setting, $V=D(\Sigma)$ and $f$ is the function given by $\deg_{\Sigma}$. We let
\begin{itemize}
  \item $\mathfrak{C}=\Amp^1 (\Sigma)$,
  \item $\Eff^1 (f) = \Eff^1 (\Sigma)$,
\end{itemize}
as introduced in Section \ref{sec lorfan}. Then $f$ is $\Amp^1 (X)$-Lorentzian, and the cones satisfy the condition required in Section \ref{pos cone lor}.

\begin{claim}
The set $\mathcal{C}=\{D_\rho| \rho\in \Sigma(1)\}$ is a generating set.
\end{claim}

By definition, $\mathcal{C}$ is a finite set and $\overline{\cone} (\mathcal{C}) = \Eff^1 (\Sigma)$. The requirement on the weak Hodge-Riemann property follows from the definition of Lorentzian fans (see Section \ref{sec lorfan}).

\subsubsection{Transcendental classes on projective manifolds}
In this setting, let $X$ be a $n$-dimensional smooth projective variety over $\mathbb{C}$.
The vector space is $V=H^{1,1} (X, \mathbb{R})$, and the function $f=\vol$ is the volume function of $(1,1)$-classes. We let
\begin{itemize}
  \item $\mathfrak{C}$ be the K\"ahler cone of $X$,
  \item $\Eff^1 (f)=\Eff^1 (X)$,
\end{itemize}
as introduced in Section \ref{sec pos}. Then $f$ is $\mathfrak{C}$-Lorentzian, and the cones satisfy the conditions required in Section \ref{pos cone lor}.

\begin{claim}
Let $\mathcal{C}=\Mov^1 (X)\cup \prim^1 (X)$ be the union of movable $(1,1)$-classes and prime divisor classes. Then $\mathcal{C}$ is a generating set.
\end{claim}

The requirement $\mathcal{C}$ generating $\Eff^1 (X)$ follows from the existence of divisorial Zariski decomposition, see Section \ref{sec zardeco}.

The finiteness of $\mathcal{C}\cap \ker \mathfrak{L}$ follows from Lemma \ref{lem mov restr}, Proposition \ref{kernel mov} and Proposition \ref{finite eff}, which extend to transcendental classes by the same arguments. The point is that we can use the result in \cite{nystromDualityMorse} as a replacement of \cite{BFJ09} to prove Proposition \ref{kernel mov} and Proposition \ref{finite eff}.

The requirement on weak Hodge-Riemann property can be checked as follows. It is clear if $D$ is given by a prime divisor, so we may assume that $D\in \Mov^1 (X)$. Since the weak Hodge-Riemann property is closed under taking limits, we may assume that $D=\pi_* A$ for some modification $\pi: \widehat{X} \rightarrow X$ and some K\"ahler class $A$ on $\widehat{X}$. Then the quadratic form can be written as
$$q(x,y) = \pi^*x\cdot \pi^*y \cdot A\cdot \pi^*\alpha_1 \cdot... \cdot \pi^*\alpha_{n-3}.$$
Define the quadratic form $\widehat{q}$ on $\widehat{X}$ by
$$\widehat{q}(x',y') = x'\cdot y' \cdot A\cdot \pi^*\alpha_1 \cdot... \cdot \pi^*\alpha_{n-3},\ x',y' \in H^{1,1}(\widehat{X}, \mathbb{R}),$$
then $\widehat{q}$ has the weak Hodge-Riemann property and $q$ is the restriction of $\widehat{q}$ to $H^{1,1}(X, \mathbb{R})$, via the injection
\begin{equation*}
  \pi^*: H^{1,1}(X, \mathbb{R})\rightarrow H^{1,1}(\widehat{X}, \mathbb{R}).
\end{equation*}
It is clear that the restriction to a subspace preserves the weak Hodge-Riemann property, therefore $q$ has the weak Hodge-Riemann property on $X$.

\subsubsection{Compact K\"ahler manifolds}
In this setting, let $X$ be a compact K\"ahler manifold of dimension $n$.
The vector space is $V=H^{1,1} (X, \mathbb{R})$, and the function $f=\vol$ is the volume function of $(1,1)$-classes. We let
\begin{itemize}
  \item $\mathfrak{C}$ be the K\"ahler cone of $X$,
  \item $\Eff^1 (f)=\Eff^1 (X)$,
\end{itemize}
as introduced in Section \ref{sec pos}. Then $f$ is $\mathfrak{C}$-Lorentzian, and the cones satisfy the conditions required in Section \ref{pos cone lor}.

\begin{claim}
Assume that every class in the collection $\mathfrak{L}=(L_1,...,L_{n-2})$ is big and nef. Let $\mathcal{C}=\Mov^1 (X)\cup \prim^1 (X)$ be the union of movable $(1,1)$-classes and prime divisor classes. Then $\mathcal{C}$ is a generating set.
\end{claim}

The first and third requirement for a generating set follow from the same argument as the transcendental case on a complex projective manifold.
The subtlety lies in the finiteness requirement in the second condition, in which we need the assumption that every class in the collection $\mathfrak{L}=(L_1,...,L_{n-2})$ is big and nef.

First we show that the analog of Proposition \ref{kernel mov} holds in the K\"ahler setting. In the proof for Proposition \ref{kernel mov}, we applied the Lefschetz hyperplane theorem. Note that on a K\"ahler manifold there may be no nontrivial subvarieties. Hence the Lefschetz hyperplane theorem is not applicable. Instead of restricting to a hyperplane, we apply the hard Lefschetz theorem for K\"ahler classes.
 We prove the result by replacing the restriction to an ample hypersurface by intersecting against a K\"ahler class.

\begin{prop}\label{kernel mov kah}
Let $X$ be a compact K\"ahler manifold of dimension $n$. Assume that $\mathfrak{L}=(L_1,...,L_{n-2})$ is a critical collection of nef classes, then
$$\Mov^1 (X)\cap \ker \mathfrak{L}=\{0\}.$$
\end{prop}

\begin{proof}
Fix a K\"ahler class $A$ on $X$, and fix $r\in [n-2]$. Let $\alpha \in \Mov^1 (X)$ satisfy 
\begin{equation}\label{movable a1}
  \mathfrak{L}\cdot \alpha=0.
\end{equation}

By the assumption (\ref{movable a1}), $ \mathfrak{L}_{\setminus r} \cdot A \cdot L_r\cdot \alpha=0$. Applying the weak Hodge-Riemann property of $\mathfrak{L}_{\setminus r} \cdot A$ and the criticality of $\mathfrak{L}$, we get that
\begin{equation*}
  \mathfrak{L}_{\setminus r} \cdot A \cdot \alpha^2 \leq 0.
\end{equation*}
On the other hand, by the definition of movable classes, $\mathfrak{L}_{\setminus r} \cdot A \cdot \alpha \in \Eff^1(X)^*$, which implies that
\begin{equation*}
  \mathfrak{L}_{\setminus r} \cdot A \cdot \alpha^2 \geq 0.
\end{equation*}
Then we obtain that $\mathfrak{L}_{\setminus r} \cdot A \cdot \alpha^2 = 0$. Using the facts
$$\mathfrak{L}_{\setminus r} \cdot A \cdot L_r ^2 >0,
\mathfrak{L}_{\setminus r} \cdot A \cdot L_r\cdot \alpha=0,$$
and applying Proposition \ref{proportional} implies that
\begin{equation*}
  \mathfrak{L}_{\setminus r} \cdot A \cdot \alpha =0.
\end{equation*}

Now taking $(\mathfrak{L}_{\setminus r}, A)$ as the new collection and inducting on the number of occurrences $A$, we finally obtain that $$A^{n-2} \cdot \alpha =0.$$ By the hard Lefschetz theorem, we get $\alpha =0$.

This finishes the proof.

\end{proof}

It remains to show that the set of prime divisor classes $[D]$ satisfying $\mathfrak{L}\cdot [D]=0$ is finite. Note that $\mathfrak{L}\cdot [D]=0$ implies that there is some $L_k$ such that
\begin{equation*}
  L_k ^{n-1} \cdot [D]=0.
\end{equation*}
Using \cite{tosattiNullLocus}, this shows that $D \subset \mathbb{B}_+ (L_k)$. Since every $L_k$ is big and nef, the set of such prime divisors is finite.
This completes the proof of our claim.

\begin{rmk}
Assuming Demailly's conjecture on transcendental Morse inequalities \cite{BDPP13}, we can drop the big nef assumption on $\mathfrak{L}$.
\end{rmk}

In summary, we obtain:

\begin{thrm}\label{local af summ}
Let $\mathfrak{L}=(L_1,...,L_{n-2})$ be a collection of nef classes. Then in either of the following situation:
\begin{itemize}
  \item $X$ is a smooth projective variety and the collection $\mathfrak{L}$ is critical,
  \item  $X$ is a compact K\"ahler manifold and every element of $\mathfrak{L}$ is big,
  \item  $\Sigma$ is a Lorentzian fan and the collection $\mathfrak{L}$ is critical,
\end{itemize}
taking the generating set as constructed above, the local Hodge index inequality holds.
\end{thrm}

\section{The numerical characterization}\label{sec char}

In this section, we apply the local Hodge index inequality to study the numerical characterization of HL classes. We focus on the case where $\mathfrak{L}=(L_1,...,L_{n-2})$ is a supercritical collection of nef classes in the geometric setting (i.e. projective varieties or compact K\"ahler manifolds).
Recall from the discussions in the introduction, the numerical characterization is closely related to the space
\begin{equation*}
  V_{\mathfrak{L}, \eff} = \Span_{\mathbb{R}} \{[D]: D\in \prim(X),\ \mathfrak{L} \cdot [D]=0\}.
\end{equation*}

\subsection{Toy example}
To illustrate the main ideas and explain why the local Hodge index inequality is useful, we start with a toy example.

Let $X$ be a 3-dimensional smooth projective variety over $\mathbb{C}$ (the same argument works for any algebraically closed field of characteristic $0$), and let $L\in N^1(X)$ be a big nef class such that $L\cdot [D]\neq 0$ for every $D\in\prim (X)$.
 For example, this holds if $L$ admits a smooth semipositive real $(1,1)$-form representatives that is strictly positive outside a finite union of curves.

\begin{claim}
The map $L: N^1(X) \rightarrow N^2(X)$ is an isomorphism.
\end{claim}

If $L$ is a semiample divisor, then the condition implies that the semiample fibration of $L$ is semismall. By \cite{cataldoMigHLsemismallmap}, $L$ is an HL class. For the nef case, this may be known to experts, but we do not find a reference.

\begin{proof}[Proof of the Claim]
It suffices to prove injectivity. Namely, for any $\alpha \in N^1(X)$ with $L\cdot \alpha=0$, we must show that $\alpha=0$. 
Note that for any prime divisor $D$, we have
\begin{equation*}
    L_{|D}^2>0 \text{ and } L_{|D}\cdot \alpha_{|D}=0.
\end{equation*}
By the Hodge index inequality on $D$, we then obtain $\alpha_{|D}^2\leq 0$.
Since $D\in\prim (X)$ is arbitrary, this implies $-\alpha^2 \in \Eff^1 (X)^*$. This is a toy model of the local Hodge index inequality.
Because $L$ lies in the interior of $\Eff^1(X)$ and $L\cdot (-\alpha^2)=0$, we obtain $\alpha^2=0$. Now take a smooth ample hypersurface $H$.
On $H$, we have
\begin{equation*}
    \alpha_{|H}^2=0 , L_{|H}^2>0 \text{ and } L_{|H}\cdot \alpha_{|H}=0.
\end{equation*}
In particular, this implies the equality
\begin{equation*}
    (L_{|H}\cdot \alpha_{|H})^2=(\alpha_{|H}^2)(L_{|H}^2).
\end{equation*}
By the Hodge index theorem on $H$, the classes $\alpha_{|H}$ and $L_{|H}$ are proportional. Thus $\alpha_{|H}=c L_{|H}$ for some constant $c\in \mathbb{R}$. Multiplying both sides by $L_{|H}$ implies $c=0$. Hence $\alpha_{|H}=0$. Then by the Lefschetz hyperplane theorem, we obtain $\alpha=0$. This finishes the proof of our claim.
\end{proof}

\subsection{Algebraic case}

To highlight the ideas in higher dimensions,
we first consider the case where all $L_k$ are big and nef.

\begin{thrm}\label{ker bignef}
Let $X$ be a smooth projective variety of dimension $n$, and let $\mathfrak{L}=(L_1,...,L_{n-2})$ be a collection of big nef classes on $X$. Then
\begin{equation*}
  \ker \mathfrak{L} = V_{\mathfrak{L}, \eff}.
\end{equation*}
In particular, $\mathfrak{L}$ is an HL class if and only if $\mathfrak{L} \cdot [D]\neq 0$ for any prime divisor $D$ on $X$.
\end{thrm}

\begin{proof}
Fix $r\in [n-2]$.
Let $\alpha\in  \ker \mathfrak{L}$. By Theorem \ref{local af summ}, there exists a class $\beta$ such that
\begin{itemize}
  \item $\beta - \alpha \in V_{\mathfrak{L}, \eff}$, and
  \item $-\beta^2 \cdot \mathfrak{L}_{\setminus r} \in \Eff^1 (X)^*$.
\end{itemize}

Since $\mathfrak{L} \cdot \beta =0$, we get that $$L_r \cdot \beta^2 \cdot \mathfrak{L}_{\setminus r} =0.$$
Since $L_r$ is big, this implies that
\begin{equation}\label{mov 0}
  \beta^2 \cdot \mathfrak{L}_{\setminus r}  =0.
\end{equation}

Let $A$ be a very ample line bundle on $X$, and let $H\in |A|$ be a general smooth hypersurface such that for every $i\in [n-2]$,
\begin{itemize}
  \item $H$ meets every component of $\mathbb{B}_+ (L_i)$ properly;
  \item the intersection of $H$ with each divisorial component of $\mathbb{B}_+ (L_i)$ is irreducible.
\end{itemize}
The existence of such a general hypersurface $H$ is guaranteed by Bertini's theorems.

On the hypersurface $H$, we have:
\begin{itemize}
  \item $\beta_{|H}^2 \cdot {\mathfrak{L}_{\setminus r}}_{|H}  =0$,
  \item $\beta_{|H}\cdot {L_r}_{|H}\cdot {\mathfrak{L}_{\setminus r}}_{|H}=0$,
  \item ${L_r}_{|H} ^2 \cdot {\mathfrak{L}_{\setminus r}}_{|H} >0$.
\end{itemize}
The first bullet follows from (\ref{mov 0}), and the second follows from $\mathfrak{L} \cdot \beta =0$. For the third, note that
\begin{equation*}
    {L_r}_{|H} ^2 \cdot {\mathfrak{L}_{\setminus r}}_{|H} ={L_r} ^2 \cdot {\mathfrak{L}_{\setminus r}}\cdot A.
\end{equation*}
By Theorem \ref{nonvanishing}, it therefore suffices to verify that the numerical dimension of the sum of any $k$-element subcollection of $(L_r,L_r,\mathfrak{L}_{\backslash r}, A)$ is at least $k$. This is immediate as every $L_i$ is big and nef.
The first and second statements yield
\begin{equation*}
    (\beta_{|H}\cdot {L_r}_{|H}\cdot {\mathfrak{L}_{\setminus r}}_{|H})^2=(\beta_{|H}^2 \cdot {\mathfrak{L}_{\setminus r}}_{|H})({L_r}_{|H} ^2 \cdot {\mathfrak{L}_{\setminus r}}_{|H}).
\end{equation*}
By Theorem \ref{prop geom}, 
there exists $c\in\mathbb{R}$ such that 
$${\mathfrak{L}_{\setminus r}}_{|H} \cdot \beta_{|H}=c{\mathfrak{L}_{\setminus r}}_{|H} \cdot L_{r|H}.$$ Multiplying both sides by $L_{r|H}$ shows that $c=0$. Hence,
$${\mathfrak{L}_{\setminus r}}_{|H} \cdot \beta_{|H} =0.$$

Note that the restriction of every $L_i$ to $H$ is also big and nef. By induction on $H$, for the class $\beta_{|H}\in \ker {\mathfrak{L}_{\setminus r}}_{|H}$ there exist prime divisors $D_i$ on $H$ such that
\begin{equation}\label{restric b}
  \beta_{|H} = \sum_{i} a_i [D_i],
\end{equation}
where each $D_i$ satisfies $${\mathfrak{L}_{\setminus r}}_{|H}\cdot [D_i] =0.$$ In particular, for some $k\in [n-2]\setminus r$,
\begin{equation*}
  {L_k}_{|H} ^{n-2}\cdot [D_i] =0.
\end{equation*}
Indeed, otherwise the restriction of every $L_k$ to $D_i$ would be big and nef, contradicting with the equality ${\mathfrak{L}_{\setminus r}}_{|H}\cdot [D_i] =0$.
Consequently, viewing $D_i$ as a cycle of codimension two in $X$, we have:
$${L_k}^{n-2}\cdot D_i =0  \text{ on }X ,$$ 
which, by definition, implies $D_i\subset \Null (L_k)$. 
Using the coincidence of the null locus and the non-K\"ahler locus \cite{nakamayeNullLocus, elmnp09restrictedvolume, tosattiNullLocus}, we obtain
\begin{equation}\label{inclusion di}
  D_i \subset \bigcup_{k\in [n-2]\setminus r} \mathbb{B}_+ (L_k).
\end{equation}

By the choice of $H$, each $D_i\in \prim(H)$ is of the form $D_i = E_i \cap H$ for some divisorial component $$E_i \subset \bigcup_{k\in [n-2]\setminus r} \mathbb{B}_+ (L_k). $$
Hence (\ref{restric b}) can be written as
\begin{equation*}
  ( \beta- \sum_{i} a_i [E_i])_{|H} =0.
\end{equation*}
By the Lefschetz hyperplane theorem, it follows that
\begin{equation}\label{eq beta}
  \beta= \sum_{i} a_i [E_i].
\end{equation}

It remains to show that $\mathfrak{L} \cdot [E_i]=0$. Since ${\mathfrak{L}_{\setminus r}}_{|H}\cdot [D_i] =0$ and $D_i = E_i \cap H$,
we have ${\mathfrak{L}_{\setminus r}}_{|H}\cdot [E_i]_{|H} =0$. In particular, this shows that
\begin{equation*}
  A^2\cdot \mathfrak{L}_{\setminus r}\cdot [E_i] =0,
\end{equation*}
which implies $\mathfrak{L}_{\setminus r}\cdot [E_i] =0$. Consequently, $\mathfrak{L} \cdot [E_i]=L_r \cdot \mathfrak{L}_{\setminus r}\cdot [E_i] =0$.

This completes the proof.

\end{proof}

\begin{exmple}
As typical examples, if every $L_k$ in the collection $\mathfrak{L}$ is big nef and $\codim \mathbb{B}_+ (L_k) \geq 2$, then $\mathfrak{L} \cdot [D]\neq 0$ for any prime divisor $D$ on $X$, which implies that $\mathfrak{L}=L_1 \cdot...\cdot L_{n-2}$ is an HL class. To the best of our knowledge, this provides the first series of two-dimensional HL classes arising as complete intersections of nef classes for which neither freeness nor semi-positivity is required.
\end{exmple}

To extend the characterization to more general collections, we need to study a kind of null locus of a general collection of nef classes. 
We start with the following result. We would like to thank Prof. Junyan Cao for the proof.

\begin{lem}\label{union locus}
Let $X$ be a smooth projective variety of dimension $n$ and let $d, k\in [n]$ with $d+k \leq n$. Assume that $L$ is a nef class on $X$ such that $\nd(L)\geq k+d$. Then the Zariski closure of the union of irreducible subvarieties $V$ of codimensional $d$ satisfying 
$$L^k \cdot [V]=0 $$
is proper in $X$.
\end{lem}
\begin{proof}
If $k+d=n$, then $L$ is big and nef. The union of such $V$ is contained in the null locus $\Null(L)$, which is a proper closed subset of $X$.

It remains to treat the case $k+d<n$. Fix a very ample line bundle $A$, and choose general sections $s_0,...,s_{n-k-d} \in H^0(X,A)$. Consider the rational map
\begin{align*}
  f: &X\dashrightarrow \mathbb{P}^{n-k-d},\\
& x\mapsto [s_0 (x):\cdots:s_{n-k-d} (x)].
\end{align*}
Let $\pi: \widehat{X} \rightarrow X$ be a resolution of $f$, and denote by $\widehat{f}: \widehat{X} \rightarrow \mathbb{P}^{n-k-d}$ the induced map. Then $\pi^* L + \widehat{f}^* c_1(\mathcal{O}(1))$ is big and nef on $\widehat{X}$. We claim that 
\begin{equation*}
    \bigcup_{L^k \cdot [V]=0,\codim V=d} V  \subset \pi \big(\Null (\pi^* L + \widehat{f}^* c_1(\mathcal{O}(1)))\cup \exc(\pi)\big).
\end{equation*}
This implies the Lemma since the right-hand side is a proper closed subset of $X$.

To prove the claim, let $V\subset X$ be a subvariety of codimension $d$ such that $L^k\cdot [V]=0$. We may assume $V\not\subset \pi(\exc(\pi))$, since otherwise the inclusion is  clear. Let $\widehat{V}$ be the proper transform of $V$ in $\widehat{X}$. Then 
$$(\pi^* L + \widehat{f}^* c_1(\mathcal{O}(1)))^{n-d}\cdot [\widehat{V}] =0$$
because $\nd_{\widehat{V}}(\pi^*L)\leq k-1$ and $\nd_{\widehat{V}}(\widehat{f}^* c_1(\mathcal{O}(1)))\leq \nd(\widehat{f}^* c_1(\mathcal{O}(1)))=n-k-d$. This implies $$\widehat{V}\subset \Null (\pi^* L + \widehat{f}^* c_1(\mathcal{O}(1))),$$ and therefore $$V\subset \pi \big( \Null (\pi^* L + \widehat{f}^* c_1(\mathcal{O}(1)))\big),$$ which completes the proof.
\end{proof}

As a corollary, we obtain:

\begin{lem}\label{locus}
Let $X$ be a smooth projective variety of dimension $n$ and let $\mathfrak{L}=(L_1,\cdots,L_{n-2})$ be a collection of nef classes.
\begin{enumerate}
  \item If the collection $\mathfrak{L}$ is critical,
  \begin{equation*}
  N_1 (\mathfrak{L})= \overline{\bigcup_{\mathfrak{L}\cdot [D]=0, \codim D=1}D}^{\zar}
  \end{equation*}
is a proper Zariski closed set of $X$, where $D$ ranges over all such prime divisors.
  \item If the collection $\mathfrak{L}$ is supercritical,
  \begin{equation*}
  N_2 (\mathfrak{L})= \overline{\bigcup_{\mathfrak{L}\cdot [V]=0, \codim V\leq 2} V}^{\zar}
  \end{equation*}
is a proper Zariski closed set of $X$, where $V$ ranges over all such prime divisors and irreducible subvarieties of codimension 2.
\end{enumerate}
\end{lem}

\begin{proof}
We first consider the critical case.
It can be derived from Lemma \ref{union locus} as follows. By Theorem \ref{nonvanishing}, a prime divisor $D$ satisfying $\mathfrak{L}\cdot [D]=0$ is equivalent to the existence of some $I\subset [n-2]$ such that $$L_I ^{|I|}\cdot [D]=0.$$
On the other hand, the criticality yields that $\nd(L_I) \geq |I|+1$. Then applying Lemma \ref{union locus} proves the result.

The proof of the supercritical case is similar. By (1), we only need to deal with those irreducible subvarieties of codimension 2.  By Theorem \ref{nonvanishing} again, a 2-codimensional subvariety $V$ satisfying $\mathfrak{L}\cdot [V]=0$ is equivalent to the existence of some $I\subset [n-2]$ such that $$L_I ^{|I|}\cdot [V]=0.$$ On the other hand, the supercriticality yields that $\nd(L_I) \geq |I|+2$. Applying Lemma \ref{union locus} then proves the result.

This finishes the proof.
\end{proof}

\begin{rmk}
Lemma \ref{locus} (1) can be also deduced directly from Proposition \ref{finite eff}.
\end{rmk}

\begin{rmk}\label{locus free}
When $\mathfrak{L}=(L_1,\cdots,L_{n-2})$ is a collection of free line bundles, by some easy geometric arguments, the properness of $N_1 (\mathfrak{L})$ and  $N_2 (\mathfrak{L})$ also extends to a compact K\"ahler manifold.

We let $\phi_I:X \rightarrow \mathbb{P}(H^0(X,L_I))$ be the Kodaira map of $L_I$ and let $Y_I=\phi_I(X)$.

We first consider the critical case.
The collection $\mathfrak{L}$ being critical implies that a general fiber of $\phi_I$ has dimension at most $n-|I|-1$, therefore
  \begin{equation*}
    Z_1 = \bigcup_{I\subset[n-2]}\phi_I^{-1}(Y_I^{\geq n-|I|})
  \end{equation*}
is a proper closed subset, where
  \begin{equation*}
    Y_I^{\geq n-|I|}=\{y\in Y_I: \dim \phi_I^{-1}(y)\geq n-|I|\}.
  \end{equation*}

Hence it suffices to show $D\subset Z_1$
for those prime divisors $D$ with $\mathfrak{L}\cdot [D] =0$.
By Theorem \ref{nonvanishing}, a prime divisor $D$ satisfies $\mathfrak{L}\cdot [D]=0$ if and only if there exists some $ I \subset [n-2]$ such that $$L_I ^{|I|} \cdot [D]=0.$$ This is equivalent to that $\dim \phi_I(D) <|I|$,
which implies that
  \begin{equation*}
    D \subset \phi_I^{-1}(Y_I^{\geq n-|I|}).
  \end{equation*}
In fact, otherwise $D \cap \phi_I^{-1}(Y_I^{\leq n-1-|I|})\neq \emptyset$, then there exists $y \in \phi_I(D)$ such that $$\dim {\phi_I}_{|D}^{-1}(y) = \dim D\cap \phi_I^{-1}(y)\leq n-1-|I|  .$$
  This is a contradiction, since by $\dim \phi_I(V) <|I|$, a general fiber of ${\phi_I}_{|D}$ has dimension at least $n-|I|$. This completes the proof of the critical case.

The supercritical case is similar. We only need to note that the supercriticality implies that
  \begin{equation*}
    Z_2 = \bigcup_{I\subset[n-2]}\phi_I^{-1}(Y_I^{\geq n-1-|I|})
  \end{equation*}
is a proper closed subset, and by similar arguments an irreducible subvariety $V$ of codimension two satisfying $\mathfrak{L}\cdot [V]=0$ implies that for some $I\subset [n-2]$,
  \begin{equation*}
    V \subset \phi_I^{-1}(Y_I^{\geq n-1-|I|}).
  \end{equation*}  
This proves the supercritical case.
\end{rmk}

 We now extend Theorem \ref{ker bignef} to the case where the collection is given a rearrangement of supercriticality. This condition is slightly stronger than supercriticality in the sense made precise in Proposition \ref{flag equiv} below.

\begin{thrm}\label{ker supercritical}
  Let $X$ be a smooth projective variety of dimension $n$, and let $\mathfrak{L}=(L_1,...,L_{n-2})$ be a collection of nef classes.
  If $\nd(L_i)\geq i+2$ for every $i$, then
  \begin{equation*}
    \ker \mathfrak{L} = V_{\mathfrak{L}, \eff}.
  \end{equation*}
\end{thrm}

\begin{proof}
Let $\alpha \in \ker (\mathfrak{L})$. By Theorem \ref{local af summ}, there exists $\beta$ such that
  \begin{itemize}
    \item $\beta-\alpha\in V_{\mathfrak{L}, \eff}$;
    \item $-\beta^2 \cdot \mathfrak{L}_{\backslash n-2} \in \Eff^1(X)^{*}$.
  \end{itemize}
  Since $L_{n-2}$ is big and $$(-\beta^2 \cdot \mathfrak{L}_{\backslash n-2})\cdot L_{n-2}=0,$$ it follows that $$-\beta^2 \cdot \mathfrak{L}_{\backslash n-2}=0.$$

By Lemma \ref{locus}, we may choose a very ample hypersurface $H$ such that
  \begin{itemize}
    \item the intersection $H\cap D$ is irreducible for any prime divisor $D\subset N_2(\mathfrak{L})$;
    \item the intersection $H\cap V$ is proper for any irreducible component $V$ of $N_2 (\mathfrak{L})$.
  \end{itemize}
  On $H$, we have:
  \begin{itemize}
    \item $\beta_{|H}^2 \cdot {\mathfrak{L}_{\backslash n-2}}_{|H}=0$;
    \item $\beta_{|H} \cdot {L_{n-2}}_{|H}\cdot {\mathfrak{L}_{\backslash n-2}}_{|H}=0$;
    \item ${L_{n-2}}_{|H}^2 \cdot {\mathfrak{L}_{\backslash n-2}}_{|H}>0$.
  \end{itemize}
  The first and the second statements are clear. The third follows from Theorem \ref{nonvanishing} -- it suffices to verify that the numerical dimension of the sum of any $k$-element subcollection of $(L_{n-2},L_{n-2},\mathfrak{L}_{\backslash n-2}, A)$ is at least $k$, which is immediate as the collection is supercritical.
Similarly as in the proof of Theorem \ref{ker bignef}, using Theorem \ref{prop geom} shows that
  \begin{equation*}
    \beta_{|H}\cdot {\mathfrak{L}_{\backslash n-2}}_{|H}=0.
  \end{equation*}

  It is clear that ${\mathfrak{L}_{\backslash n-2}}_{|H}=({L_1}_{|H},...,{L_{n-3}}_{|H})$ satisfies
  \begin{equation*}
    \nd({L_i}_{|H})\geq i+2,\ \forall i \in [n-3].
  \end{equation*}
  So by induction on $H$, $\beta_{|H}$ can be written as
  \begin{equation*}
    \beta_{|H}=\sum a_i D_i
  \end{equation*}
  for some $a_i \in \mathbb{R}$ and prime divisors $D_i$ in $H$ such that ${\mathfrak{L}_{\backslash n-2}}_{|H}\cdot D_i=0$. In particular, viewed as a codimension-two subvariety of $X$, we have
\begin{equation*}
  {\mathfrak{L}_{\backslash n-2}}\cdot D_i =0
\end{equation*}

 Hence $D_i \subset N_2(\mathfrak{L})$ by the definition of $N_2(\mathfrak{L})$.
 Now by the choice of $H$, $D_i$ cannot be an irreducible component of $N_2(\mathfrak{L})$.
  Hence there exists a prime divisor $E_i\subset N_2(\mathfrak{L})$ containing $D_i$. Note that $D_i\subset H\cap E_i$ and $H\cap E_i$ being irreducible imply $D_i=H\cap E_i$.
  Thus we have $\beta_{|H}=\sum a_i {E_i}_{|H}$. It follows that
  \begin{equation*}
    \beta=\sum a_i E_i
  \end{equation*}
  by Lefschetz hyperplane theorem. We also have $\mathfrak{L}\cdot E_i =0$, since
  \begin{equation*}
    H\cdot \mathfrak{L}\cdot E_i={\mathfrak{L}}_{|H}\cdot {E_i}_{|H} ={L_{n-2}}_{|H}\cdot {\mathfrak{L}_{\backslash n-2}}_{|H}\cdot D_i=0.
  \end{equation*}

  Up to now, we have shown that $\beta \in V_{\mathfrak{L}, \eff}$ and so is $\alpha$ by construction of $\beta$.
  Hence we get
  \begin{equation*}
    \ker \mathfrak{L} \subset  V_{\mathfrak{L}, \eff}.
  \end{equation*}

  The converse inclusion is clear. This completes the proof.
\end{proof}

\subsubsection{Supercritical collections under a rearrangement}
Recall that a collection $\mathcal{F}$ of subsets of $[n-2]$ is called a full (increasing) flag of $[n-2]$ if
\begin{equation*}
 I_1 \subset I_2 \subset...\subset I_{n-2} =[n-2],
\end{equation*}
and $|I_k| =k, \forall k\in[n-2]$.
Given a supercritical collection of nef classes $\mathfrak{L}=(L_1,...,L_{n-2})$, associated with a full flag $\mathcal{F}$ of $[n-2]$, one can construct a new collection of nef classes
\begin{equation*}
  \mathfrak{L}_\mathcal{F}=(L_{I_1},...,L_{I_{n-2}}).
\end{equation*}
It is clear that the collection $\mathfrak{L}_\mathcal{F}$ is also supercritical and satisfies that $\nd(L_{I_k})\geq k+2$. Then by Theorem \ref{ker supercritical}, we have
\begin{equation*}
  \ker \mathfrak{L}_\mathcal{F} =V_{\mathfrak{L}_\mathcal{F}, \eff}.
\end{equation*}

We show that the space $V_{\mathfrak{L}, \eff}$ of an arbitrary supercritical collection $\mathfrak{L}$ is given by the description of $\ker \mathfrak{L}_\mathcal{F}$ by taking into account of all the full increasing flags $\mathcal{F}$ of $[n-2]$.

\begin{prop}\label{flag equiv}
Let $\mathfrak{L}=(L_1,...,L_{n-2})$ be a supercritical collection of nef classes, then the following equality holds:
\begin{equation*}
  \sum_{\mathcal{F}} \ker \mathfrak{L}_\mathcal{F} = \sum_{\mathcal{F}} V_{\mathfrak{L}_\mathcal{F}, \eff}=V_{\mathfrak{L}, \eff},
\end{equation*}
where the sum is taken over all full flags $\mathcal{F}$ of $[n-2]$.
\end{prop}

\begin{proof}
By Theorem \ref{ker supercritical}, for any full flag $\mathcal{F}$ of $[n-2]$, $\ker \mathfrak{L}_\mathcal{F}$ is spanned by prime divisors $D$ such that $\mathfrak{L}_\mathcal{F} \cdot [D]=0$. The inclusion
\begin{equation*}
  \sum_{\mathcal{F}} \ker \mathfrak{L}_\mathcal{F} \subset V_{\mathfrak{L}, \eff}=\Span_{\mathbb{R}} (\prim(X) \cap \ker \mathfrak{L}),
\end{equation*}
follows once we verify that any such $D$ also satisfies $\mathfrak{L} \cdot [D]=0$. This holds since for any ample (or K\"ahler class) $A$,
\begin{equation*}
  0\leq \mathfrak{L} \cdot [D]\cdot A \leq \mathfrak{L}_\mathcal{F} \cdot [D]\cdot A.
\end{equation*}

For the converse direction, let $D$ be a prime divisor such that $\mathfrak{L} \cdot [D]=0$. By Theorem \ref{nonvanishing}, this is equivalent to the existence of some $I\subset [n-2]$ such that $$L_I ^ {|I|} \cdot [D]=0.$$ 
Denote $|I|=s$. Consider the full flag $\mathcal{F}_I$ containing $I$ as a piece:
\begin{equation*}
  I_1 \subset...\subset I_s =I\subset ...\subset [n-2],
\end{equation*}
then
\begin{equation*}
  \nd_D (\sum_{i=1}^s L_{I_i}) =\nd_D (L_{I_s}) < |I_s|.
\end{equation*}
By Theorem \ref{nonvanishing} again, $\mathfrak{L}_{\mathcal{F}_I} \cdot [D]=0$.
By Theorem \ref{ker supercritical}, this proves the converse inclusion.

This finishes the proof.

\end{proof}

\begin{rmk}
Combining Theorem \ref{ker supercritical}, we see that Conjecture \ref{main conj} is equivalent to the statement:
\begin{quote}
 \emph{Let $\mathfrak{L}=(L_1,...,L_{n-2})$ be a supercritical collection of nef classes, then
\begin{equation*}
  \sum_{\mathcal{F}} \ker \mathfrak{L}_\mathcal{F} = \ker \mathfrak{L},
\end{equation*}
where the sum is taken over all full flags $\mathcal{F}$ of $[n-2]$.}
\end{quote}

\end{rmk}

\subsection{K\"ahler case}

Note that in the K\"ahler case, the Lefschetz hyperplane theorem may not be available since a K\"ahler manifold may not admit any ample hypersurface. We present an approach which works for a collection of big nef classes in the K\"ahler setting. It also provides an alternative approach to Theorem \ref{ker bignef}.

\begin{thrm}\label{kernel Kahler}
  Let $X$ be a compact K\"ahler manifold of dimension $n$ and $\mathfrak{L}=(L_1,...,L_{n-2})$ be a collection of big and nef classes. Then
  \begin{equation*}
    \ker \mathfrak{L} = V_{\mathfrak{L}, \eff}.
  \end{equation*}
\end{thrm}

\begin{proof}
  Fix a K\"ahler class $\omega$ on $X$, and for $i=0,...,n-2$ denote
$$\mathfrak{L}^{(i)}=(L_1,...,L_i,\omega,...,\omega).$$
We prove by induction on the number of inserted copies of $\omega$  that
  \begin{equation*}
    \ker \mathfrak{L}^{(i)} = V_{\mathfrak{L}^{(i)}, \eff} =\Span_{\mathbb{R}} \{[D]: D\in \prim(X),\ \mathfrak{L}^{(i)} \cdot [D]=0\}
  \end{equation*}
holds for any $i=0,...,n-2$.

Note that $i=0$ is just the usual hard Lefschetz theorem for K\"ahler classes.

By induction, we may suppose that it holds for $i=k-1$.
  Let $\alpha \in \ker \mathfrak{L}^{(k)}$. By Theorem \ref{local af summ}, we can find $\beta$ such that
  \begin{itemize}
    \item $\beta-\alpha\in \Span_{\mathbb{R}} \{[D]: D\in \prim(X),\ \mathfrak{L}^{(k)} \cdot [D]=0\}$;
    \item $-\beta^2 \cdot \mathfrak{L}^{(k)}_{\backslash k} \in \Eff^1(X)^{*}$.
  \end{itemize}
  Now clearly $(-\beta^2 \cdot \mathfrak{L}^{(k)}_{\backslash k})\cdot L_{k}=0$ by the choice of $\beta$. Since $L_k$ is big, this implies 
  \begin{equation*}
    -\beta^2 \cdot \mathfrak{L}^{(k)}_{\backslash k}=0.
  \end{equation*}
  Therefore,
  \begin{itemize}
    \item $\beta^2 \cdot \mathfrak{L}^{(k)}_{\backslash k} \cdot \omega =\beta^2 \cdot \mathfrak{L}^{(k-1)}=0$;
    \item $\beta \cdot L_{k} \cdot \mathfrak{L}^{(k-1)}=0$;
    \item $L_{k}^2 \cdot \mathfrak{L}^{(k-1)}>0$.
  \end{itemize}
  By Theorem \ref{prop geom} and the same argument as in Theorem \ref{ker bignef}, this implies 
  \begin{equation*}
    \beta \cdot \mathfrak{L}^{(k-1)}=0.
  \end{equation*}

By the inductive hypothesis, this shows: $\beta\in V_{\mathfrak{L}^{(k-1)}, \eff} = \ker \mathfrak{L}^{(k-1)}$.
  We then get the result by the following easy fact:
  \begin{align*}
   V_{\mathfrak{L}^{(k-1)}, \eff} \subset V_{\mathfrak{L}^{(k)}, \eff}.
  \end{align*}

This completes the proof.

\end{proof}

\begin{rmk}
When $L_1=...=L_{n-2}$, the above result can be considered as the K\"ahler geometric analogue of \cite[Theorem 2.2]{handelMinkExtr}, which holds for convex bodies beyond polytopes. This case corresponds to the extremals of Minkowski's quadratic inequality, see also \cite[Theorem 14.6]{handel2022AFextremals} for more results.
\end{rmk}

\subsection{Bergman fans}
Theorem \ref{kernel Kahler} extends to Lorentzian fans satisfying hard Lefschetz property. In particular, this works for the Bergman fans of (poly)matroids.

\begin{thrm}\label{kernel Bergman fan}
  Let $\Sigma$ be the Bergman $n$-fan of a matroid or polymatroid and let $\mathfrak{L}=(L_1,...,L_{n-2})$ be a collection of big and nef classes, then
  \begin{equation*}
    \ker \mathfrak{L} = V_{\mathfrak{L}, \eff}.
  \end{equation*}
\end{thrm}

Here, $V_{\mathfrak{L}, \eff}$ for a Bergman fan is the linear space $\Span_{\mathbb{R}}\{[D_\rho]: \mathfrak{L}\cdot D_\rho =0, \rho \in \Sigma(1)\}$.
The proof is identical to that of Theorem \ref{kernel Kahler}, replacing the hard Lefschetz theorems and Hodge–Riemann relations for K\"ahler manifolds by those for Bergman fans. We omit the details.

\subsection{Application}\label{sec hodge extre}

Let $\mathfrak{L}=(L_1,...,L_{n-2})$ be a collection of nef classes and let $A, B$ be nef classes. By the Hodge index inequality (also known as the Khovanski-Teissier inequality),
\begin{equation*}
  \left(A\cdot B\cdot \mathfrak{L}\right)^2 \geq \left(A^2\cdot \mathfrak{L}\right) \left(B^2\cdot \mathfrak{L}\right).
\end{equation*}

It is natural to ask when equality holds, i.e., to characterize the extremals of this inequality.
If $A\cdot B\cdot \mathfrak{L} =0$, then the equality holds automatically since the right-hand side is nonnegative. By Theorem \ref{nonvanishing} (or Theorem \ref{PosiCri} in a more general setting), this situation is completely characterized in terms of numerical dimensions.
Hence it suffices to consider the case $A\cdot B\cdot \mathfrak{L} >0$. By Theorem \ref{prop geom}, the equality holds if and only if $$\mathfrak{L} \cdot (A-cB)=0$$ for some $c>0$. Therefore, the problem of describing the extremals reduces to understanding the kernel $\ker \mathfrak{L}$.

Under the same assumptions, the Hodge index inequality also implies that the sequence$$\{a_k=A^k \cdot B^{n-k}\}_{0\leq k\leq n}$$ is log-concave, that is,
\begin{equation*}
  a_k ^2 \geq a_{k-1} a_{k+1}, \ \forall 1\leq k \leq n-1.
\end{equation*}

Many important log-concave sequences arising in combinatorics and geometry can be obtained in this way \cite{huhICM2018}. It is an important and interesting question to characterize the extremals of $ a_k ^2 = a_{k-1} a_{k+1}$ for a fixed $k$,
see \cite{chan2021logconcave, handel2022AFextremals} for recent developments.

In the algebro-geometric setting, when $X$ is a smooth projective variety (or a compact K\"ahler manifold), it follows from \cite{BFJ09} (or \cite{fuxiaoTeissProp}) that:
\begin{thrm}
If $A, B$ are big and nef on $X$, then $$a_{n-1} ^n = a_n ^n a_0$$ if and only if $A, B$ are proportional.
\end{thrm}
Note that $a_{n-1} ^n = a_n ^n a_0$ is equivalent to the condition $a_k ^2 =a_{k-1} a_{k+1}$ for every $ 1\leq k \leq n-1$.

It remains open\footnote{This question had also been asked to the second named author by Professor J.-P. Demailly in 2014 at Institut Fourier.} to characterize $a_l ^2 =a_{l-1} a_{l+1}$ for a single index $l$. Our main results provide a complete answer to this question for $A, B$ big and nef.

\begin{thrm}
Let $X$ be a smooth projective variety of dimension $n$, and let $\mathfrak{L}=(L_1,...,L_{n-2})$ be a collection of nef classes such that $\nd(L_i)\geq i+2$ for every $i\in[n-2]$. Let $A, B$ be nef classes such that
\begin{equation*}
  0< \left(A\cdot B\cdot \mathfrak{L}\right)^2 = \left(A^2\cdot \mathfrak{L}\right) \left(B^2\cdot \mathfrak{L}\right),
\end{equation*}
then $$A-cB =\sum_{i} c_i [D_i]$$
for some prime divisors $D_i$ satisfying $\mathfrak{L}\cdot [D_i]=0$. Here $c>0$ is the constant such that $\mathfrak{L} \cdot A=c\mathfrak{L} \cdot B$.

In particular, if $A, B$ are big and nef, then
\begin{equation*}
  (A^k \cdot B^{n-k})^2 = (A^{k-1} \cdot B^{n-k+1})(A^{k+1} \cdot B^{n-k-1})
\end{equation*}
holds for some $k\in[n-1]$ if and only if
\begin{equation*}
  A-cB =\sum_{i} c_i [D_i]
\end{equation*}
where $D_i$ are prime divisors satisfying $$A^{k-1} \cdot B^{n-k-1}\cdot [D_i]=0$$
and $c>0$ is given by $A^{k} \cdot B^{n-k-1}=cA^{k-1} \cdot B^{n-k}$.
\end{thrm}

\begin{proof}
This follows straightforward from Theorem \ref{ker supercritical}.
\end{proof}

\begin{rmk}
This result can be seen as the algebraic analogue of \cite[Theorem 14.6]{handel2022AFextremals} (see also \cite[Theorem 2.2]{handelMinkExtr}).
By Theorems \ref{kernel Kahler}, \ref{kernel Bergman fan}, the above characterization of extremals for big and nef $A, B$ also extends to the K\"ahler setting and to Lorentzian fans satisfying the hard Lefschetz theorems and Hodge-Riemann relations.
\end{rmk}

\section{Examples and further discussions}\label{sec arbi}

When the collection $\mathfrak{L}$ is arbitrary, the characterization of $\ker \mathfrak{L}$ is much more subtle and new inputs are required to attack the problem. Nevertheless,
in this section we present some positive evidence and examples and illustrate why we formulate Conjecture \ref{main conj1} as stated in the introduction. We also give some brief discussions on what could be expected beyond Hodge index theory.

Recall the conjectural picture for a general collection:

\begin{conj}\label{conjec sec7}
Let $X$ be a smooth projective variety (or compact K\"ahler manifold) of dimension $n$. Let $\mathfrak{L}=(L_1,...,L_{n-2})$ be a collection of nef classes on $X$. Define the vector spaces $V_{\mathfrak{L},\eff}$ and $V_{\mathfrak{L},\deg}$ as follows:
\begin{align*}
V_{\mathfrak{L}, \eff} &= \Span_{\mathbb{R}} \{[D]: D\in \prim(X),\ \mathfrak{L} \cdot [D]=0\},\\
V_{\mathfrak{L},\deg}&={\Span}_{\mathbb{R}}\{\mu_*(\widetilde{\alpha}-\widetilde{\beta}): (\widetilde{\alpha}, \widetilde{\beta})\ \text{is a $\mu^*\mathfrak{L}$-degenerate pair on $\widetilde{X}$}\},
\end{align*}
where the linear span for $V_{\mathfrak{L},\deg}$ is taken over all birational modifications $\mu : \widetilde{X}\rightarrow X$ and all $\mu^*\mathfrak{L}$-degenerate pairs $(\widetilde{\alpha}, \widetilde{\beta})$ on $\widetilde{X}$.
Then $$\ker \mathfrak{L} = V_{\mathfrak{L},\eff}+V_{\mathfrak{L},\deg}.$$
In particular, the class $\mathfrak{L}$ is an HL class if and only if $V_{\mathfrak{L},\deg}=\{0\}$ and $\mathfrak{L} \cdot [D] \neq 0$ for any prime divisor on $X$.
\end{conj}

As noted in the introduction,
when the collection $\mathfrak{L}$ is not subcritical, by Theorem \ref{nonvanishing} or \ref{PosiCri} the complete intersection class $$\mathfrak{L}=L_1\cdot...\cdot L_{n-2}= 0,$$
thus $\ker \mathfrak{L}$ is the whole space and Conjecture \ref{conjec sec7} holds trivially in this case, since any element of $N^1 (X)$ can be written as a difference of two ample classes. Therefore, Conjecture \ref{conjec sec7} is nontrivial only when the collection is subcritical.

\begin{rmk}
Assume that the collection is subcritical. Let $M$ be the unique maximal subcritical subset of $\mathfrak{L}$ (see Corollary \ref{max subcritical}). Assume that
\begin{enumerate}
  \item \emph{the class $\mathfrak{L}_M$ can be represented by a (smooth) irreducible subvariety $V$};
  \item \emph{the restriction to $V$ induces an isomorphism from $N^1 (X)$ to $N^1 (V)$.}
\end{enumerate}
By Proposition \ref{subcr to cri}, the restricted collection $\mathfrak{L}_{\setminus M |V}$ is critical on $V$. Assume that $\alpha\in \ker\mathfrak{L}$, then by (2) we have $\alpha_{|V} \in \ker\mathfrak{L}_{\setminus M |V}$. Therefore, with the above assumptions the characterization can be reduced to the critical case on a lower dimensional variety.
\end{rmk}


Simple examples show that the degenerate pair is an indispensable obstruction for HL classes.

\begin{exmple}

Let $X=\mathbb{P}^1 \times \mathbb{P}^1 \times \mathbb{P}^1$ and let $L=\pi_1 ^*\mathcal{O}_{\mathbb{P}^1}(1)+\pi_2 ^*\mathcal{O}_{\mathbb{P}^1}(1)$ where $\pi_k$ is the projection to the $k$-th factor. It is easy to check that $L: N^1 (X)\rightarrow N^2 (X)$ is not an isomorphism. In this case, $\nd(L)=2$ and there is no prime divisor $D$ such that $L\cdot [D]=0$, and
\begin{align*}
&L\cdot \pi_1 ^*\mathcal{O}_{\mathbb{P}^1}(1)\cdot \pi_2 ^*\mathcal{O}_{\mathbb{P}^1}(1)=0,\\
&L\cdot (\pi_1 ^*\mathcal{O}_{\mathbb{P}^1}(1)- \pi_2 ^*\mathcal{O}_{\mathbb{P}^1}(1))=0.
\end{align*}
Hence, $(\pi_1 ^*\mathcal{O}_{\mathbb{P}^1}(1), \pi_2 ^*\mathcal{O}_{\mathbb{P}^1}(1))$ is a degenerate pair for $L$.

\end{exmple}

\begin{prop}\label{vdeg}
Let $\mathfrak{L}=(L_1,...,L_{n-2})$ be a collection of nef classes on $X$, then $$V_{\mathfrak{L}, \deg}\subset \ker \mathfrak{L}.$$
Assume further that the collection $\mathfrak{L}$ is supercritical, then $V_{\mathfrak{L}, \deg}=\{0\}$.
\end{prop}

\begin{proof}

Assume that $(\widetilde{\alpha}, \widetilde{\beta})$ is a $\mu^*\mathfrak{L}$-degenerate pair on the birational modification $\mu: \widetilde{X}\rightarrow X$, that is, for some ample class $A$ on $\widetilde{X}$ we have that
\begin{equation*}
  \mu^*\mathfrak{L}\cdot \widetilde{\alpha}\cdot \widetilde{\beta}=0,\
\mu^*\mathfrak{L}\cdot \widetilde{\alpha}\cdot A=\mu^*\mathfrak{L}\cdot \widetilde{\beta}\cdot A.
\end{equation*}

We first consider the first statement.
If $\mu^*\mathfrak{L}\cdot \widetilde{\alpha}=0$, then the latter equality shows that $\mu^*\mathfrak{L}\cdot \widetilde{\beta}=0$, which implies
$\mu^*\mathfrak{L}\cdot (\widetilde{\alpha}-\widetilde{\beta})=0$.
Next we assume that both $\mu^*\mathfrak{L}\cdot \widetilde{\alpha}$ and $\mu^*\mathfrak{L}\cdot \widetilde{\beta}$ are nonzero, by Theorem \ref{prop geom} the equality $\mu^*\mathfrak{L}\cdot \widetilde{\alpha}\cdot \widetilde{\beta}=0$ implies that $\mu^*\mathfrak{L}\cdot \widetilde{\alpha}=c \mu^*\mathfrak{L}\cdot \widetilde{\beta}$ for some $c>0$. The requirement on the intersection numbers against $A$ forces $c=1$. We also get $\mu^*\mathfrak{L}\cdot (\widetilde{\alpha}-\widetilde{\beta})=0$.

In either case, we have that
\begin{equation*}
  \mu^*\mathfrak{L}\cdot (\widetilde{\alpha}-\widetilde{\beta})=0.
\end{equation*}
By the projection formula, this shows that $\mu_*  (\widetilde{\alpha}-\widetilde{\beta}) \in \ker \mathfrak{L}$, which proves that $V_{\mathfrak{L}, \deg}\subset \ker \mathfrak{L}.$

Next we consider the case when the collection $\mathfrak{L}$ is supercritical, then the collection $\mu^*\mathfrak{L}$ is supercritical for any birational modification $\mu$. By Theorem \ref{nonvanishing} the condition $\mu^*\mathfrak{L}\cdot \widetilde{\alpha}\cdot \widetilde{\beta}=0$ ensures that either $(\widetilde{\alpha}+\widetilde{\beta})^2 =0$ or at least one of $\widetilde{\alpha}, \widetilde{\beta}$ is zero. In the latter case, if $\widetilde{\alpha}=0$, then $\mu^*\mathfrak{L} \cdot \widetilde{\beta} =0$ and the supercriticality shows that $\widetilde{\beta}=0$. Hence we may assume $\widetilde{\alpha}, \widetilde{\beta}$ are nonzero and
$(\widetilde{\alpha}+\widetilde{\beta})^2 =0$, then $$\widetilde{\alpha}\cdot \widetilde{\beta} =\widetilde{\alpha}^2 =\widetilde{\beta}^2 =0.$$
By Theorem \ref{prop geom}, this shows that $\widetilde{\alpha} = a \widetilde{\beta}$ for some $a>0$. Using $\mu^*\mathfrak{L} \cdot A \cdot \widetilde{\alpha}=\mu^*\mathfrak{L} \cdot A \cdot \widetilde{\beta}$ yields that $\widetilde{\alpha}-\widetilde{\beta}=0$. Thus, $V_{\mathfrak{L}, \deg}=\{0\}$.

This finishes the proof.
\end{proof}

Next we show that the statement of Conjecture \ref{conjec sec7} has birational invariance.

\begin{prop}
Let $\pi: \widehat{X}\rightarrow X$ be a birational morphism. Assume that Conjecture \ref{conjec sec7} holds on $\widehat{X}$, then it also holds on $X$.
\end{prop}

\begin{proof}
Assume that $\alpha \in \ker \mathfrak{L}$. Then $\pi^*\alpha \in \ker \pi^*\mathfrak{L}$. By the assumption on $\widehat{X}$, there exist prime divisors $\widehat{D}_j$ on $\widehat{X}$ with $\pi ^* \mathfrak{L} \cdot [\widehat{D}_j] =0$  and $( \pi\circ \mu_k) ^* \mathfrak{L}$-degenerate pairs $(\alpha_k, \beta_k)$ on the modification $\mu_k: X_k \rightarrow \widehat{X}$ such that
\begin{equation*}
  \pi^* \alpha = \sum_{j} a_j [\widehat{D}_j] + \sum_k b_k \mu_{k*} (\alpha_k - \beta_k).
\end{equation*}
Applying the pushforward $\pi_*$, we get that
\begin{equation*}
 \alpha = \sum_{j} a_j [\pi_* \widehat{D}_j] + \sum_k b_k (\pi \circ \mu_{k})_* (\alpha_k - \beta_k),
\end{equation*}
which clearly implies that $\alpha \in V_{\mathfrak{L}, \eff}+V_{\mathfrak{L},\deg}$ on $X$.

The argument for the supercritical collection is similar, since the pullback of a supercritical collection via a birational morphism remains supercritical.

\end{proof}

\subsection{Toric varieties}
We briefly explain (under certain assumptions) why Conjecture \ref{conjec sec7} should hold on a toric variety using the result of Shenfeld and van Handel.
This had already been noted in \cite[Section 16]{handel2022AFextremals}.
We provide a few more details.

Recall that Shenfeld and van Handel proved the following remarkable result \cite[Theorem 2.13]{handel2022AFextremals}:

\begin{thrm}\label{AFextremals}
  Let $\mathcal{C}=(C_1,...,C_{n-2})$ be a collection of polytopes in $\mathbb{R}^n$, and let $A,B$ be convex bodies such that $V_n(A,B,\mathcal{C})>0$. Then
  \begin{equation*}
    V_n(A,B,\mathcal{C})^2=V_n(A,A,\mathcal{C})V_n(B,B,\mathcal{C})
  \end{equation*}
  if and only if there exists $a>0,v\in \mathbb{R}^n$, and a finite number of $\mathcal{C}$-degenerate pairs $$(M_1,N_1),...,(M_m,N_m),$$ so that $A+M_1+\cdots+M_m$ and $aB+v+N_1+\cdots+N_m$ have the same supporting hyperplane in all $(\mathbf{B},\mathcal{C})$-extreme normal directions, or equivalently,
  \begin{equation*}
    h_{A+M_1+\cdots +M_m}=h_{aB+v+N_1+\cdots+N_m}\ \text{a.e.-} S_{\mathbf{B},\mathcal{C}}.
  \end{equation*}
Here, $h_K \in C^0 (S^{n-1})$ denotes the support function of the convex body $K$.
\end{thrm}

Here, we use the following notations:
\begin{itemize}
  \item $\mathbf{B}$ is the unit ball in $\mathbb{R}^n$.
  \item $V_n (-)$ denotes the mixed volume of convex bodies.
  \item a pair $(M,N)$ of convex bodies is called $\mathcal{C}$-degenerate if
  \begin{equation*}
    V_n(M,N,\mathcal{C})=0, V_n(M,\mathbf{B},\mathcal{C})=V_n(N,\mathbf{B},\mathcal{C}).
  \end{equation*}
  \item $S_{\mathbf{B},\mathcal{C}}$ is the mixed area measure.
  \item Since $\mathcal{C}$ consists of polytopes, the set of all $(\mathbf{B},\mathcal{C})$-extreme normal directions is given by
  \begin{equation*}
    \supp S_{\mathbf{B},\mathcal{C}}=\{u\in S^{n-1}: \dim F(C_I,u)\geq |I|, \forall I\subset [n-2]\} ,
  \end{equation*}
  where $F(C_I,u)$ is the face of $C_I=\sum_{i\in I}C_i$ with outer normal direction $u$.
\end{itemize}

We briefly recall some basics on toric varieties
(see e.g. \cite{fultonToricBook, coxToricBOOK}) and use the following notations:
\begin{itemize}
  \item $\Sigma$ is a projective simplicial fan in $\mathbb{R}^n$. Equivalently, $\Sigma$ is the outer normal fan of a rational simple polytope $P$.
  \item $X=X_{\Sigma}$ is the $\mathbb{Q}$-smooth toric variety associated to the fan $\Sigma$.
  \item $A^*(X)$ is the Chow ring of $X$ with coefficients in $\mathbb{R}$. Note that for complete simplicial toric varieties, we have
        \begin{equation*}
          A^*(X)\cong N^*(X) \cong H^*(X,\mathbb{R}).
        \end{equation*}
        $A^1(X)$ can be interpreted explicitly as follows:
        \begin{equation*}
          A^1(X)\cong \frac{PL( \Sigma)}{L(\Sigma)}.
        \end{equation*}

  \item $\{D_{\rho}:\rho\in \Sigma(1)\}$ is the set of all $T$-invariant prime divisors. It is also a generating set of the pseudoeffective cone of $X$, i.e.,
        \begin{equation*}
          \Eff^1(X)=\overline{\cone}\{[D_{\rho}]: \rho\in \Sigma(1)\}.
        \end{equation*}
        As a conewise linear function, $D_{\rho}$ is defined by
        \begin{equation*}
          u_{\rho}\mapsto 1, u_{\sigma} \mapsto 0, \text{for those } \sigma\neq \rho.
        \end{equation*}
        Here $u_{\rho}$ is the primitive generator of the ray $\rho$.
  \item Let $Q$ be a polytope  homothetic to a Minkowski summand of $P$. Its support function $h_Q$ is conewise linear and convex on $\Sigma$. So it determines a nef class $[Q]:=[h_Q]$ in $A^1(X)$. In fact, all nef classes take this form, i.e.,
        \begin{equation*}
          \Nef^1(X)=\{[Q]: \exists \lambda \in \mathbb{R}_{> 0 } \ s.t. \ P= \lambda Q+Q_1\}.
        \end{equation*}
  \item $\mathfrak{L}=(L_1,...,L_{n-2})$ is a collection of $T$-invariant nef divisors on $X$.
  \item $\mathcal{C}=(C_1,...,C_{n-2})$ is the collection of polytopes associated to $\mathfrak{L}$.
\end{itemize}

We take $\alpha \in \ker \mathfrak{L}$ a rational class. Note that $\alpha$ can be written as a difference of two rational ample classes $\beta$ and $\gamma$. Thus we have
  \begin{equation*}
    (\beta\cdot \gamma \cdot \mathfrak{L})^2=(\beta\cdot \beta \cdot \mathfrak{L})(\gamma\cdot \gamma \cdot \mathfrak{L}).
  \end{equation*}
  Let $A,B$ be the polytopes associated to $\beta,\gamma$. By the  Bernstein-Kushnirenko-Khovanskii theorem, the above equality is equivalent to
  \begin{equation*}
    V_n(A, B ,\mathcal{C})^2=V_n(A, A ,\mathcal{C})V_n(B, B ,\mathcal{C}).
  \end{equation*}
Now by Theorem \ref{AFextremals}, there exist $v\in \mathbb{R}^n$ and $\mathcal{C}$-degenerate pairs $(M_1,N_1),...,(M_m,N_m)$ such that
  \begin{equation}\label{sptfcn}
    h_{A+M_1+\cdots +M_m}=h_{B+v+N_1+\cdots+N_m}\ \text{a.e.-} S_{\mathbf{B},\mathcal{C}}.
  \end{equation}
Here the constant $a=1$ since $\mathfrak{L}\cdot \beta=\mathfrak{L}\cdot \gamma$.

\emph{Assume that the degenerate pairs $(M_i,N_i)$ can be taken to be homothetic to rational polytopes}.
We may choose a projective simplicial fan $\widetilde{\Sigma}$, which is a subdivision of $\Sigma$, such that $h_{M_i},h_{N_i}$ are all conewise linear on $\widetilde{\Sigma}$ \cite{coxToricBOOK}.
  Denote by $$\mu:\widetilde{X}=X_{\widetilde{\Sigma}} \rightarrow X$$ the toric modification determined by the subdivision $\widetilde{\Sigma} \rightarrow \Sigma$.
  Since $h_{A+M_1+\cdots +M_m}$ and $h_{B+v+N_1+\cdots+N_m}$ are conewise linear on $\widetilde{\Sigma}$, they are determined by their values on $\{u_{\rho}\}_{\rho \in \widetilde{\Sigma}(1)}$.
  Then by (\ref{sptfcn}) and the definition of $D_{\rho}$, we can find $a_{\rho}$ such that
  \begin{equation}\label{sptrelation}
    h_{A+M_1+... +M_m}-h_{B+v+N_1+...+N_m}= \sum_{\rho\in \widetilde{\Sigma}(1) \backslash \supp S_{\mathbf{B},\mathcal{C}}} a_{\rho}D_{\rho}.
  \end{equation}

We use $[-]_{\widetilde{X}}, [-]_X$ to denote classes on $\widetilde{X}$ and $X$ respectively.
By the additivity of the map $K \mapsto h_K$, by (\ref{sptrelation}) we get that
  \begin{equation*}
    ([A]_{\widetilde{X}}-[B]_{\widetilde{X}}) +([M_1]_{\widetilde{X}}-[N_1]_{\widetilde{X}})+...+ ([M_m]_{\widetilde{X}} -[N_m]_{\widetilde{X}} )
   \end{equation*}
is an element of the space $\Span_{\mathbb{R}}\{ [D_u]_{\widetilde{X}}:u\in \widetilde{\Sigma}(1)\backslash \supp S_{\mathbf{B},\mathcal{C}}\}.$
Since $$[A]_{\widetilde{X}}=\mu^*[A]_{X}, [B]_{\widetilde{X}}=\mu^*[B]_{X},$$ we find that
  \begin{equation*}
    \alpha = [A]_{X}-[B]_{X}=\mu_*([A]_{\widetilde{X}}-[B]_{\widetilde{X}}).
  \end{equation*}
By the definition of $M_i, N_i$, it is clear that $\mu_*([M_i]_{\widetilde{X}}-[N_i]_{\widetilde{X}}) \in V_{\mathfrak{L}, \deg}$.
  Now the result follows by the following observation:
  \begin{equation*}
    \mu_*\Span_{\mathbb{R}}\{ [D_u]_{\widetilde{X}}:u\in \widetilde{\Sigma}(1)\backslash \supp S_{\mathbf{B},\mathcal{C}}\}=V_{\mathfrak{L},\eff}.
  \end{equation*}

\subsection{Smooth fibration satisfying Leray-Hirsch}

We consider a special critical case given by a smooth fibration of relative dimension one.
Let $\pi: X\rightarrow  Y$ be a smooth fibration between compact K\"ahler manifolds with $\dim X =n$ and $\dim Y =n-1$. Assume that the Leray-Hirsch theorem holds for this fibration, in particular, there exist classes $$\xi \in H^{1,1}(X, \mathbb{R}),\ g_1,...,g_t \in H^{1,0}(X)$$
such that on every fiber $F$ their restrictions form a basis of $H^{1,1}(F, \mathbb{R})$ (resp. $H^{1,0}(F)$), and
$H^{1,1}(X, \mathbb{R})$ is generated by $\pi^* H^{1,1}(Y, \mathbb{R}), \xi$ and $H^{1,0}(X) \cdot \pi^* H^{0,1}(Y)$.

Let $h_1,...,h_s \in H^{1,0}(Y)$ be a basis on the base manifold.

Consider the collection $\mathfrak{L}=(\pi^* \omega_1,...,\pi^* \omega_{n-2})$ where $\omega_i$'s are K\"ahler classes on $Y$. This is a critical collection $\mathfrak{L}$ since $\nd(\pi^* \omega_i)=n-1$ for every $i$. We show that Conjecture \ref{main conj1} holds in this setting. 

\begin{prop}
    In the setting as above, we have
    $\ker\mathfrak{L}= V_{\mathfrak{L}, \deg}$.
\end{prop}
\begin{proof}
    We claim that $\ker \mathfrak{L}\subset \pi^*H^{1,1}(Y,\mathbb{R})$, which immediately implies the identity $\ker\mathfrak{L}= V_{\mathfrak{L}, \deg}$. To see this, let $B\in \ker \mathfrak{L}$ and write $B=\pi^* \gamma$ for some $\gamma \in H^{1,1}(Y, \mathbb{R})$. Decompose $\gamma$ into the difference of K\"ahler classes $\gamma_1-\gamma_2$,  we find that $B$ is the difference of an $\mathfrak{L}$-degenerate pair $(\pi^*\gamma_1,\pi^*\gamma_2).$

    It remains to prove the claim.
    By the Leray-Hirsch assumption, we can decompose $B$ as follows:
\begin{equation*}
  B=\pi^* b + c \xi + \sum_{i\in [t], j\in [s]} a_{ij} g_i \cdot \pi^* \overline{h_j} +\overline{\sum_{i\in [t], j\in [s]} {a_{ij}} { g_i} \cdot \pi^*\overline{ h_j}}.
\end{equation*}
The condition $\mathfrak{L}\cdot B =0$ implies that, for any $\alpha\in H^{1,1}(Y, \mathbb{R})$, $$\mathfrak{L}\cdot B \cdot \pi^* \alpha =0.$$ Computing the left-hand side gives
\begin{equation*}
  \mathfrak{L}\cdot B \cdot \pi^* \alpha = c (\omega_1\cdot...\cdot \omega_{n-1}\cdot \alpha) (\pi_* \xi),
\end{equation*}
so $c=0$. 
Next, for any $p, q$, intersecting with $g_p \cdot \pi^* \overline{h_q}$ gives 
\begin{equation*}
  \mathfrak{L}\cdot B \cdot (g_p \cdot\pi^* \overline{h_q}) = \sum_{i\in [t], j\in [s]}  \pi_* (g_p \cdot \overline{g_i}) \overline{a_{ij}} (\mathfrak{L}\cdot \pi^*h_j \cdot \overline{\pi^*h_q}) =0.
\end{equation*}
Define the matrices $$G=[\pi_* (g_p \cdot \overline{g_i})]_{p,i}, H=[\mathfrak{L}\cdot \pi^*h_j \cdot \overline{\pi^*h_q}]_{j, q}.$$
By the Hodge-Riemann relations on the fiber (a curve) and on $Y$, the matrices $G, H$ are definite. Hence the coefficient matrix $A=[a_{ij}]$ must vanish. Thus $B=\pi^*b$, proving the claim.
\end{proof}

\subsection{Compact complex torus}

We explain why Conjecture \ref{main conj1} holds on a compact complex torus, essentially following Panov's result \cite{Panov1987ONSP}.

Let $X=\mathbb{C}^n / \Gamma$ be a compact complex torus of dimension $n$. Then 
$$H^{1,1} (X, \mathbb{R})\cong \Lambda^{1,1}_{\mathbb{R}}(\mathbb{C}^n),$$
the space of real $(1,1)$-forms with constant coefficients on $\mathbb{C}^n$. Under this identification, both the pseudo-effective cone and nef cone of $(1,1)$-classes coincide with the cone of semi-positive real $(1,1)$-forms.

\begin{prop}
    For any subcritical collection $\mathfrak{L}$ on $X$, $\ker \mathfrak{L}=V_{\mathfrak{L}, \deg}$.
\end{prop}

\begin{proof}
When $\mathfrak{L}$ is critical, the result reduces essentially to that of Panov \cite{Panov1987ONSP}.

It remains to address the subcritical but non-critical case. In this setting, there exists a non-empty subset $I\subset [n-2]$ such that
$\nd(L_I)=|I|$. By the results in Section \ref{sec collec}, we may assume without loss of generality that $I$ is maximal (i.e. there is no $I\subsetneq J \subset [n-2]$ such that $\nd(L_J)=|J|$) and, by relabeling, $I=[k]$.
Since $\nd(L_I)=|I|$, we may choose a coordinates system $(z_1,...,z_n)$ for $\mathbb{C}^n$ such that $$L_I=\sum_{i=1}^{k}\sqrt{-1}dz_i\wedge d\overline{z_i}.$$
Let $W=\{z_1=...=z_{k}=0\}$ be the subspace of $\mathbb{C}^n$.
One easily verifies the decomposition $$\ker\mathfrak{L}=\ker \mathfrak{L}_{I^c|W}+ \ker(\Lambda^{1,1}_{\mathbb{R}}(\mathbb{C}^n) \xrightarrow{restriction} \Lambda^{1,1}_{\mathbb{R}}(W) ),$$
where $$\mathfrak{L}_{I^c|W}=(L_{k+1|W},...,L_{n-2|W})$$ is a collection of $n-k-2$ semi-positive real $(1,1)$-forms on the $(n-k)$-dimensional space $W$. Here, we view $\ker \mathfrak{L}_{I^c|W}\subset \Lambda^{1,1}_{\mathbb{R}}(W)$ as a subspace of $\Lambda^{1,1}_{\mathbb{R}}(\mathbb{C}^n)$ in a natural way -- $\Lambda^{1,1}_{\mathbb{R}}(W)$ consists of those $(1,1)$-forms involving no $dz_i$ or $d\bar{z}_j$ for $1\leq i, j\leq k$.
By the maximality of $I$ and Proposition \ref{subcr to cri}, the collection $\mathfrak{L}_{I^c|W}$ is critical. Panov's result then implies $$\ker \mathfrak{L}_{I^c|W} = V'_{\deg},$$
where $V'_{\deg}$ is the space generated by the differences of all $\mathfrak{L}_{I^c|W}$-degenerate pairs. Since $V'_{\deg}\subset V_{\mathfrak{L}, \deg}$ is immediate, it suffices to show that
$$\ker(\Lambda^{1,1}_{\mathbb{R}}(\mathbb{C}^n) \xrightarrow{restriction} \Lambda^{1,1}_{\mathbb{R}}(W) ) \subset V_{\mathfrak{L},\deg}.$$
To this end, note that $\ker(\Lambda^{1,1}_{\mathbb{R}}(\mathbb{C}^n) \xrightarrow{restriction} \Lambda^{1,1}_{\mathbb{R}}(W) )$ is generated by forms of the type
\begin{equation*}
 \alpha= \sum_{i,j\in[k]}a_{i,j}\sqrt{-1}dz_i\wedge d\overline{z_j}+b \sqrt{-1}dz_l\wedge d\overline{z_s} +\overline{b} \sqrt{-1} dz_s\wedge d\overline{z_l},
\end{equation*}
where $A=[a_{ij}]$ is a Hermitian matrix, $b\in \mathbb{C}$, $l\in [k] $ and $s\in[n]\backslash [k]$.
We define
\begin{align*}
  &\beta=\alpha+c(\sum_{i=1}^{k}\sqrt{-1}dz_i\wedge d\overline{z_i}+\sqrt{-1}dz_s\wedge d\overline{z_s}), \\
  &\gamma=c(\sum_{i=1}^{k}\sqrt{-1}dz_i\wedge d\overline{z_i}+\sqrt{-1}dz_s\wedge d\overline{z_s}),
\end{align*}
where $c>0$ a constant chosen large enough to ensure that $\beta$ is semi-positive.
Since $\nd(L_I+\beta+\gamma) \leq |I|+1$, Theorem \ref{nonvanishing} implies that $\mathfrak{L}_I\cdot \beta \cdot \gamma=0.$
Consequently, $(\beta,\gamma)$ is a $\mathfrak{L}$-degenerate pair. As $\alpha=\beta-\gamma$, we conclude that $\alpha \in V_{\mathfrak{L},\deg}$, which completes the proof.
\end{proof}

\subsection{Beyond Hodge index}
The main results of this paper focus on Hodge index theory and log-concavity, i.e. HL classes of dimension two. As this approach does not  extend trivially to classes of other dimensions, several interesting questions arise. Motivated by \cite{hushangxiao2023hardlefArxiv} and the results established herein, we provide a brief discussion on potential generalizations.

By introducing ``$m$-lefness'', a  notion of partial positivity for algebraic maps defined via the defect of semismallness, 
we proved the following result in our previous work:

\begin{thrm}[see \cite{hushangxiao2023hardlefArxiv}]
Let $X$ be a smooth projective variety of dimension $n$, and let $\mathfrak{L}=(L_1, ...,L_{n-p-q})$ be a collection of free line bundles on $X$. 
If $L_I$ is $(|I|+p+q)$-lef for any $I\subset [n-p-q]$, then $\mathfrak{L}$ is an HL class on $H^{p, q}(X)$.
\end{thrm}

Recall that a free line bundle $L$ is $m$-lef on $X$ if and only if, for any irreducible subvariety $T\subset X$,
\begin{equation*}
  2\dim T -\dim \Phi_L (T) -n \leq n-m,
\end{equation*}
where $\Phi_L: X \rightarrow \mathbb{P}(H^0 (X, L))$ is the Kodaira map associated with $L$.
This condition can be expressed purely numerically. Specifically, $L_I$ being $|I|+p+q$-lef is equivalent to the following: for any irreducible subvariety $V$,
\begin{equation*}
 \dim \Phi_{L_I} (V)\geq (|I|+p+q)-2\codim V.
\end{equation*}
In particular, whenever $(|I|+p+q)-2\codim V \geq 0$, we have
\begin{equation*}
  L_I ^{(|I|+p+q)-2\codim V} \cdot [V] \neq 0.
\end{equation*}

Now, consider the case when $p=q$. For any irreducible subvariety $V$ of codimension $p$, the condition simplifies to 
$$L_I ^{|I|} \cdot [V] \neq 0.$$
By Theorem \ref{nonvanishing}, the condition $L_I ^{|I|} \cdot [V] \neq 0$ for all $I\subset [n-2p]$ is equivalent to 
\begin{equation*}
\mathfrak{L}\cdot [V]=L_1\cdot...\cdot L_{n-2p}\cdot [V] \neq 0.
\end{equation*}
Since we only consider subvarieties of codimension $p$, this represents only a subset of the requirements for $L_I$ being $|I|+2p$-lef. However, geometric intuition from the $p=1$ case suggests this non-vanishing may be the primary obstruction to the HL property on $N^{p} (X)$ or $H^{p,p}(X)$. Furthermore, for a compact complex torus of dimension $n$ and a collection of nef classes $\mathfrak{L}=(L_1, ...,L_{n-p-q})$, it was shown in \cite{huxiaohardlef2022Arxiv} that $\mathfrak{L}$ is an HL class on $H^{p,q} (X)$ if and only if $$\nd(L_I)\geq |I|+p+q$$ for all $I\subset [n-p-q]$. 

These observations lead to the following expectation:

\begin{expec}\label{ques1}
Let $X$ be a smooth projective variety of dimension $n$, and let $\mathfrak{L}=(L_1, ...,L_{n-2p})$ be a collection of nef classes on $X$ such that $$\nd(L_I)\geq |I|+2p$$ for any $I\subset [n-2p]$.
Then the kernel of $\mathfrak{L}=L_1\cdot...\cdot L_{n-2p}$, acting as a map on $N^p (X)$ (or $H^{p,p} (X)$), is spanned by those cycle classes of irreducible subvarieties $V$ of codimension $p$ such that
\begin{equation*}
\mathfrak{L}\cdot [V] =L_1\cdot...\cdot L_{n-2p}\cdot [V]= 0.
\end{equation*}
In particular, $\mathfrak{L}$ is an HL class on $N^p (X)$ (or $H^{p,p} (X)$) if and only if
$\mathfrak{L}\cdot [V]\neq 0$
for every irreducible subvariety $V$ of codimension $p$.
\end{expec}

This characterization is highly non-trivial, as it generalizes Grothendieck's standard conjectures regarding the hard Lefschetz theorem for algebraic cycles to the setting where $\mathfrak{L}$ is a collection of nef (rather than purely ample) classes.

\section*{Declaration} 

\textbf{Data Availability:} Not applicable.

\textbf{Conflict of interest:} On behalf of all authors, the corresponding author states that there is no conflict of interest.

\bibliography{lefcharreference}
\bibliographystyle{amsalpha}

\bigskip

\bigskip

\noindent
\textsc{Tsinghua University, Beijing 100084, China}\\
\noindent
\verb"Email: hujj22@mails.tsinghua.edu.cn"\\
\noindent
\verb"Email: jianxiao@tsinghua.edu.cn"

\end{document}